\documentclass{amsart}
\usepackage{graphicx}
\usepackage[all]{xy}
\newtheorem{theorem}{Theorem}[section]

\newtheorem{lemma}[theorem]{Lemma}

\theoremstyle{definition}
\newtheorem{definition}[theorem]{Definition}
\newtheorem{Notation}[theorem]{Notation}
\theoremstyle{remark}
\newtheorem{remark}[theorem]{Remark}

\numberwithin{equation}{section}

\begin{document}

\title[On the Finite-Dimensional Irreducible Representations of $PSL_2(\mathbb{Z})$]
{On the Finite-Dimensional Irreducible Representations of
$PSL_2(\mathbb{Z})$}
\author{Melinda G. Moran and Matthew J. Thibault}
\maketitle
\begin{abstract}
We classify up to equivalence all finite-dimensional irreducible
representations of $PSL_2(\mathbb{Z})$ whose restriction to the
commutator subgroup is diagonalizable.
\end{abstract}
\section{Introduction}
In this paper, we study finite-dimensional representations of the projective
modular group $PSL_2(\mathbb{Z})$.  Adriaenssens and Le Bruyn [1]
have recently analyzed specific families of representations of
$PSL_2(\mathbb{Z})$, from the point of view of noncommutative algebraic geometry.
Also, a complete classification of finite-dimensional irreducible representations of dimension
$\leq 5$ follows from the work of Tuba and Wenzl [4] on representations of $B_3$. Our
main result completely classifies up to equivalence all
finite-dimensional irreducible representations of
$PSL_2(\mathbb{Z})$ whose restrictions to the commutator subgroup are
diagonalizable. These representations are of dimension $1$, $2$,
$3$, and $6$. \\
\indent Recall that $PSL_2(\mathbb{Z})$ has the presentation $\langle x,y$ $|$ $x^2 = y^3 = 1 \rangle$.
Its commutator subgroup has index 6 and is generated by the elements $xyxy^2, xy^2xy$;
see [3]. Since the index is 6, it follows from standard Clifford Theory that the
dimensions of the irreducible representations we are studying divide 6;
see, for example, [2, 2.7] for further explanation.  This fact is used in our analysis. \\
\indent The inspiration of this paper comes from the theory of
highest weight modules. The finite-dimensional irreducible representations of a complex semisimple Lie
algebra are diagonalizable over its Cartan subalgebra.  In our study below,
we view the commutator of $PSL_2(\mathbb{Z})$ as playing a role analogous
to that of the Cartan subalgebra. \\
\indent Another way to view our work is as follows: Set $G'$ as the commutator subgroup of $G$
and $G''$ as the second commutator subgroup of G.  It is well known
(see [3]) that $G'$ is a free group in the two generators $xyxy^2$ and $xy^2xy$, and so $G'/G''$ is a free
abelian group on two generators.  In particular, the irreducible representations of $G'/G''$ are all one-dimensional.
Therefore, by Clifford theory, the restrictions to $G'/G''$ of irreducible finite-dimensional representations of $G/G''$
must be diagonalizable.  Conversely, every finite-dimensional irreducible representation of $G/G''$ lifts to a
representation of $G$ whose restriction to $G'$ naturally factors through $G/G''$.  We see that our results
amount to a complete classification, up to equivalence, of the irreducible finite-dimensional representations
of $G/G''$.  Moreover, since $G/G''$ is abelian-by-finite, the irreducible representations are all finite-dimensional. \\
\indent Our approach is mostly elementary, relying on basic linear algebraic computations
 and case-by-case analyses. Toward the end of the paper,
extensive Maple calculations are used to determine the
6-dimensional representations; an appendix of the Maple code used is included. \\
\indent  Our main results are stated in section 4. Sections 5 and 6 are devoted to the proofs.
Preliminary results and notation are given in sections 2 and 3. \\ \\
\noindent{\bf Acknowledgements:} This work was begun at the Temple University 2004
Research Experience for Undergraduates, funded by NSF REU Site Grant DMS-0138991.
We would like to thank our REU supervisor Professor E.
Letzter for formulating the question, for his advice on this
research project, and for his help in writing this paper.  We are also happy to acknowledge helpful discussion
of this work with D. Vogan.
\section{Definitions and Notation}
We begin with some relevant definitions.  Let $\textit{k}$ denote an
algebraically closed field and $M_n(\textit{k})$ denote the set of
all $n \times n$ matrices with entries in $\textit{k}$.  Let
$GL_n(\textit{k})$ denote the set of all invertible elements of
$M_n(\textit{k})$.
\begin{definition}
We will say the ordered m-tuple $(A_1, A_2, \ldots, A_m)$, $A_i \in
M_n(\textit{k})$ for $i = 1, 2, \ldots, m$, is \textit{irreducible}
if every element of $M_n(\textit{k})$ can be written as a
$\textit{k}$-linear combination of products in the $A_i$'s, $i = 1,
2, \ldots m$. \\Also, if $X_i, X'_i \in M_n(\textit{k})$ for $i = 1,
2, \ldots m$, then $(X'_1, X'_2, \ldots, X'_m)$ is
\textit{equivalent} to $(X_1, X_2, \ldots X_m)$, denoted $(X'_1,
X'_2, \ldots, X'_m) \approx (X_1, X_2, \ldots X_m)$, if there exists
$Q \in GL_n(\textit{k})$ such that $X'_i = QX_iQ^{-1}$ for all $i =
1, 2, \ldots, m$.
\end{definition}

\begin{Notation}
\noindent The following notation will remain in effect for the entire paper:
 Let $G = PSL_2(\mathbb{Z})$, which we identify with $\langle x,y$ $|$ $x^2 =
y^3 = 1 \rangle$. Let $\rho : G \rightarrow GL_n(\textit{k})$ be an
irreducible representation of $G$.  Set $X = \rho(x)$, $Y =
\rho(y)$, $\Lambda = XYXY^2$, and $\Gamma = XY^2XY$. It then follows
that $X^2 = Y^3 = I$, and $(X,Y)$ is irreducible. Denote the entry
in the $i^{th}$ row and $j^{th}$ column of a matrix $X$ by
$X_{i,j}$. We will use $\langle f_1, \ldots, f_t \rangle$ to denote the ideal (in a given ring)
generated by  $f_1, f_2, \ldots f_t$.
\end{Notation}

\section{Preliminary Results}
\begin{remark}
Since $\textit{k}$ is algebraically closed, all irreducible
solutions to $XY=YX$ in $M_n(\textit{k})$ are one-dimensional.
\end{remark}
\begin{lemma}
All irreducible representations of $PSL_2(\mathbb{Z})$ with $\Lambda
= \Gamma$, or $\Lambda$ or $\Gamma$ a scalar matrix are
one-dimensional.
\end{lemma}
\begin{proof}
If $\Lambda = \Gamma$, then $XY = Y^2X \Lambda Y^2 = Y^2X \Gamma Y^2
= YX$.  For an arbitrary constant $c$, computation with $\Lambda =
cI$ yields
$$XY = (Y^2X)^2(XY)\Lambda Y^2 =(Y^2X)^2 \Lambda(XY)Y^2 = YX.$$
Similarly, computation with $\Gamma = cI$ yields $$XY =
(Y^2X)(XY)\Gamma (Y^2XY^2) = (Y^2X)\Gamma (XY)(Y^2XY^2)= YX.$$
Remark 3.1 thus concludes this proof.
\end{proof}

\begin{remark}
Assume $\Lambda$ and $\Gamma$ are $n \times n$ diagonal matrices, where $\Lambda_{i,i} =: \lambda_i$
and $\Gamma_{i,i} =: \gamma_i$ for $i = 1, 2, \ldots, n$.  We
observe the following properties of $\Lambda$ and $\Gamma$:
\begin{itemize}
\begin{item}
$\Lambda X \Lambda = X$, $\Gamma X \Gamma = X$, $\Lambda Y \Gamma
= Y$, $\Gamma Y^2 \Lambda = Y^2$
\end{item} \\
\begin{item}
$\Lambda \Gamma = \Gamma \Lambda$, so $(XY)^6 = (YX)^6 = I$.
\end{item} \\
\begin{item}
Since $Y$ has at least one nonzero entry per row, $\Lambda$ and
$\Gamma$ are conjugate, and $\Lambda Y \Gamma = Y$, then
$\frac{1}{\lambda_i} = \gamma_j = \lambda_k$ for each $i = 1, 2,
\ldots, n$,  and some $j$, $k \in \{1, 2, \ldots, n\}$. Note that
$\Lambda$, $\Gamma \in GL_n(\textit{k})$, so $\lambda_i, \gamma_i
\neq 0$ for each $i = 1, 2, \ldots, n$.
\end{item}
\end{itemize}
\end{remark}
The proof of the following lemma is routine and omitted.
\begin{lemma}
 Let the following properties hold:
 \\\begin{enumerate}
\begin{item}$X,Y \in GL_n(\textit{k})$ and
have exactly one nonzero entry per row and column
\end{item}
\\ \begin{item}  $(X,Y)$ is irreducible \end{item}
\\ \begin{item}  $X^2=Y^3=I$  \end{item}
\\ \begin{item}  $\Lambda$, $\Gamma$ are diagonal matrices
\end{item}
\end{enumerate}
Then $(X,Y)$ satisfies the above properties if and only if
$(PXP^{-1},PYP^{-1})$ satisfies the above properties, where $P$ is a
non-singular weighted permutation matrix (i.e. $P$ has exactly one
nonzero entry per row and column).
\end{lemma}

\section{Main Results}
Let $\rho$ be an irreducible representation of $PSL_2(\mathbb{Z})$
which maps the commutator subgroup of $PSL_2(\mathbb{Z})$ to
diagonal matrices in $GL_n(\textit{k})$. Since the index of the commutator subgroup
is 6, it follows from standard Clifford Theory (see, e.g., [2,2.7] that the
dimension of $\rho$ divides 6. We thus analyze the cases when $n = 1, 2, 3$, and $6$.
\begin{theorem}
Let $\textit{k}$ be an algebraically closed field, and let $\zeta$ be
a primitive cube root of unity if $\textit{k}$ is not of
characteristic 3.  Let $\rho: PSL_2(\mathbb{Z}) \rightarrow
GL_n(\textit{k})$ be an irreducible representation of
$PSL_2(\mathbb{Z}) = \langle x,y$ $|$ $x^2 = y^3 = 1 \rangle$ which
maps the commutator subgroup of $PSL_2(\mathbb{Z})$ to diagonal
matrices in $GL_n(\textit{k})$.
\\
i) If $\textit{k}$ is not of characteristic 2 or 3, then
$(\rho(x),\rho(y))$ is equivalent to one of the following:
\begin{enumerate}
\begin{item}
$(1,1), (-1,1), (1,\zeta), (-1,\zeta), (1,\zeta^2), (-1,\zeta^2)$
\end{item}
\\
\begin{item}
$(\left( \begin{matrix} 0 & 1 \\ 1 & 0 \end{matrix} \right),
\left( \begin{matrix} 1 & 0 \\ 0 & \zeta \end{matrix} \right))$,
$(\left( \begin{matrix} 0 & 1 \\ 1 & 0 \end{matrix} \right),
\left( \begin{matrix} 1 & 0 \\ 0 & \zeta^2 \end{matrix} \right))$,
$(\left( \begin{matrix} 0 & 1 \\ 1 & 0 \end{matrix} \right),
\left( \begin{matrix} \zeta & 0 \\ 0 & \zeta^2 \end{matrix}
\right))$
\end{item}
\\
\begin{item}
$\Bigl(\pm \left( \begin{matrix} 1 & 0 & 0 \\ 0 & -1 & 0 \\ 0 & 0 &
-1
\end{matrix} \right), \left( \begin{matrix} 0 & 1 & 0 \\ 0 & 0 & 1 \\ 1 & 0 & 0
\end{matrix} \right) \Bigr)$
\end{item}
\\
\begin{item}
$\Biggl( \left(
\begin{matrix} 0 & 1 & 0 & 0 & 0 & 0 \\ 1 & 0 & 0 & 0 & 0 & 0 \\ 0
& 0 & 0 & 1 & 0 & 0 \\ 0 & 0 & 1 & 0 & 0 & 0 \\ 0 & 0 & 0 & 0 & 0
& 1 \\ 0 & 0 & 0 & 0 & 1 & 0 \end{matrix} \right), \left(
\begin{matrix} 0 & 0 & 0 & 0 & c_1 & 0 \\ 0 & 0 & 0 & 0 & 0 & 1 \\
c_2 & 0 & 0 & 0 & 0 & 0 \\ 0 & 1 & 0 & 0 & 0 & 0 \\ 0 & 0 &
\frac{1}{c_1c_2} & 0 & 0 & 0 \\ 0 & 0 & 0 & 1 & 0 & 0 \end{matrix}
\right) \Biggr)$, \\ for $c_1,c_2 \in$ $\textit{k}$, $c_1, c_2 \neq 0$,
$(c_1,c_2) \neq (1,1)$, $(-1,1)$, $(1,-1)$, $(-1,-1)$,
$(\zeta,\zeta)$, $(\zeta^2,\zeta^2)$.
\end{item}
\end{enumerate}
ii) If $\textit{k}$ is of characteristic 2, then $(\rho(x),\rho(y))$
is equivalent to one of the following:
\begin{enumerate}
\begin{item}
$(1,1), (1,\zeta), (1,\zeta^2)$
\end{item}
\\
\begin{item}
$(\left( \begin{matrix} 0 & 1 \\ 1 & 0 \end{matrix} \right),
\left( \begin{matrix} 1 & 0 \\ 0 & \zeta \end{matrix} \right)),
(\left( \begin{matrix} 0 & 1 \\ 1 & 0 \end{matrix} \right), \left(
\begin{matrix} 1 & 0 \\ 0 & \zeta^2 \end{matrix} \right)), (\left(
\begin{matrix} 0 & 1 \\ 1 & 0 \end{matrix} \right), \left(
\begin{matrix} \zeta & 0 \\ 0 & \zeta^2 \end{matrix} \right))$
\end{item}
\\
\begin{item}
$\Biggl( \left(
\begin{matrix} 0 & 1 & 0 & 0 & 0 & 0 \\ 1 & 0 & 0 & 0 & 0 & 0 \\ 0
& 0 & 0 & 1 & 0 & 0 \\ 0 & 0 & 1 & 0 & 0 & 0 \\ 0 & 0 & 0 & 0 & 0
& 1 \\ 0 & 0 & 0 & 0 & 1 & 0 \end{matrix} \right), \left(
\begin{matrix} 0 & 0 & 0 & 0 & c_1 & 0 \\ 0 & 0 & 0 & 0 & 0 & 1 \\
c_2 & 0 & 0 & 0 & 0 & 0 \\ 0 & 1 & 0 & 0 & 0 & 0 \\ 0 & 0 &
\frac{1}{c_1c_2} & 0 & 0 & 0 \\ 0 & 0 & 0 & 1 & 0 & 0 \end{matrix}
\right) \Biggr)$, \\ for $c_1,c_2 \in$ $\textit{k}$, $c_1, c_2 \neq 0$,
$(c_1,c_2) \neq (1,1)$, $(\zeta,\zeta)$, $(\zeta^2,\zeta^2)$.
\end{item}
\end{enumerate}
iii) If $\textit{k}$ is of characteristic 3, then
$(\rho(x),\rho(y))$ is equivalent to one of the following:
\begin{enumerate}
\begin{item}
$(1,1), (-1,1)$
\end{item}
\\
\begin{item}
$\Bigl(\pm \left( \begin{matrix} 1 & 0 & 0 \\ 0 & -1 & 0 \\ 0 & 0 &
-1
\end{matrix} \right), \left( \begin{matrix} 0 & 1 & 0 \\ 0 & 0 & 1 \\ 1 & 0 & 0
\end{matrix} \right) \Bigr)$
\end{item}
\\
\begin{item}
$\Biggl( \left(
\begin{matrix} 0 & 1 & 0 & 0 & 0 & 0 \\ 1 & 0 & 0 & 0 & 0 & 0 \\ 0
& 0 & 0 & 1 & 0 & 0 \\ 0 & 0 & 1 & 0 & 0 & 0 \\ 0 & 0 & 0 & 0 & 0
& 1 \\ 0 & 0 & 0 & 0 & 1 & 0 \end{matrix} \right), \left(
\begin{matrix} 0 & 0 & 0 & 0 & c_1 & 0 \\ 0 & 0 & 0 & 0 & 0 & 1 \\
c_2 & 0 & 0 & 0 & 0 & 0 \\ 0 & 1 & 0 & 0 & 0 & 0 \\ 0 & 0 &
\frac{1}{c_1c_2} & 0 & 0 & 0 \\ 0 & 0 & 0 & 1 & 0 & 0 \end{matrix}
\right) \Biggr)$, \\ for $c_1,c_2 \in $ $\textit{k}$, $c_1, c_2 \neq 0$,
$(c_1,c_2) \neq (1,1)$, $(-1,1)$, $(1,-1)$, $(-1,-1)$.
\end{item}
\end{enumerate}
Furthermore, for $n<6$, each of the cases i), ii), and iii) yields
distinct equivalence classes. For $n=6$ and $c_1, c_2$ satisfying
the above criteria, the following correspond to the same equivalence
class: $(c_1,c_2)$, $(\frac{1}{c_1}$,$\frac{1}{c_2})$,
$(c_2,\frac{1}{c_1c_2})$, $(\frac{1}{c_2},c_1c_2)$,
$(\frac{1}{c_1c_2},c_1)$, and $(c_1c_2,\frac{1}{c_1})$. If
$(c'_1,c'_2)$ does not equal any of the preceding pairs, then
$(c'_1,c'_2)$ represents an equivalence class distinct from
$(c_1,c_2)$.
\end{theorem}
\noindent The remainder of the paper is devoted to the proof of this theorem.
\begin{remark}
\indent Since the commutator subgroup of $G$ is generated by
$xyxy^2$ and $xy^2xy$, $\rho$ maps the commutator subgroup of $G$ to
diagonal matrices in $GL_n(\textit{k})$ if and only if $\Lambda$ and
$\Gamma$ are diagonal matrices.  Thus the problem reduces to finding
distinct equivalence classes of $(X,Y)$ where $X := \rho(x)$, $Y :=
\rho(y)$, $X^2 = I = Y^3$, $\Lambda = XYXY^2$ is a diagonal matrix,
$\Gamma = XY^2XY$ is a diagonal matrix, and $(X,Y)$ is irreducible.
\end{remark}
\section{Cases when n=1, n=2, and n=3}
In this section, we prove (i), (ii), and (iii) of Theorem 4.1,
considering separately the cases when $n=1, 2,$ and $3$.
\subsection{n=1}
For $n = 1$, $\Lambda$ and $\Gamma$ are trivially diagonal.
We only require that $(X,Y) =
(a,b)$ where $a^2 = 1 = b^3$. \\
\indent Distinct $(X,Y)$ create distinct equivalence classes since
$(X',Y') \approx (X,Y)$ requires that $(X',Y') =$
$(QXQ^{-1},QYQ^{-1})=(X,Y)$. In fields of characteristic 2, the only
solution to $X^2 = 1$ is $X=1$. In fields of characteristic other
than 2, $X^2 = 1$ has the 2 distinct solutions $X = \pm 1$. In
fields of characteristic 3, the only solution to $Y^3 = 1$ is $Y=1$.
In fields of characteristic other than 3, $Y^3 = 1$ has the 3
distinct solutions $Y = 1, \zeta, \zeta^2$. We therefore arrive at
the desired result.
\subsection{n=2}
 Let $\Lambda =
\left( \begin{matrix} \lambda_1 & \\ & \lambda_2
\end{matrix} \right)$. \smallskip First assume $\lambda_1 = \lambda_2$. Then $\Lambda = \lambda_1 I$,
 which is a contradiction by Lemma 3.2.  Thus assume $\lambda_1 \neq \lambda_2$.
By explicitly solving $\Lambda X \Lambda = X$, we find that for any
$i, j \in \{1,2\}$, $\lambda_i \lambda_j = 1$ if $X_{i,j} \neq 0$.
If $X$ has more than one nonzero entry per row or column, it follows
that $\lambda_1 = \lambda_2$, which is a contradiction to our
assumption. Therefore since $X$ is non-singular, $X$ has exactly one
nonzero entry per row and column.
\\
\indent Since $\Lambda$ and $\Gamma$ are conjugate, we find $\Gamma$
has diagonal entries $\gamma_1$, $\gamma_2$ where $\gamma_1 \neq
\gamma_2$.  Solving $\Lambda Y \Gamma = Y$, we find that for any $i,
j \in \{1,2\}$, $\lambda_i \gamma_j = 1$ if $Y_{i,j} \neq 0$.  If
$Y$ has more than one nonzero entry per row or column, then either
$\lambda_1 = \lambda_2$ or $\gamma_1 = \gamma_2$, both of which are
contradictions. Therefore since $Y$ is non-singular, $Y$
has exactly one nonzero entry per row and column. \smallskip Hence $X$ and $Y$
are of the form $\left(
\begin{matrix} * & 0 \\ 0 & * \end{matrix} \right)$ or $\left(
\begin{matrix} 0 & * \\ * & 0 \end{matrix} \right)$. If $Y$ is of
the latter type, then $Y^3 \neq I$, so $Y$ must be diagonal. Clearly $X$
is of the latter form by the irreducibility of $(X,Y)$. Using $X^2 =
I$, we see that $X = \left(
\begin{matrix} 0 & x_1 \\ \frac{1}{x_1} & 0
\end{matrix} \right)$ for some nonzero $x_1 \in \textit{k}$.
Conjugating $X$ and $Y$ by the weighted permutation matrix $P = \left(
\begin{matrix} 1 & 0 \\ 0 & x_1 \end{matrix} \right)$, we see by
Lemma 3.4 that $(X,Y) \approx (X_1,Y_1) := (\left( \begin{matrix} 0 & 1 \\
1 & 0
\end{matrix} \right),\left( \begin{matrix} a & 0 \\ 0 & b
\end{matrix} \right))$, where $(X_1,Y_1)$ satisfies the hypotheses of Lemma 3.4.
By Remark 3.1, $a \neq b$. Also, if $a \neq b$,
\begin{align*}
\frac{1}{a-b}(X_1Y_1-bX_1) & = \left(\begin{matrix} 0 & 0 \\
1 & 0
\end{matrix} \right) & \frac{1}{a-b}(aX_1^2-Y_1) & = \left(\begin{matrix} 0 & 0 \\ 0 & 1
\end{matrix} \right) \\
\frac{1}{a-b}(aX_1-X_1Y_1) & = \left(\begin{matrix} 0 & 1 \\ 0 & 0
\end{matrix} \right) & \frac{1}{a-b}(Y_1-bX_1^2) & = \left(\begin{matrix} 1 & 0 \\ 0 & 0
\end{matrix} \right).
\end{align*}
Thus $(X_1,Y_1)$ is irreducible if and only if $a \neq b$.
\\ \indent Finally, $a^3 = b^3 = 1$ by $Y_1^3 = I$. In fields of characteristic 3, $a = b = 1$ since
$a^3 = b^3 = 1$, which is a contradiction. Thus there are no
irreducible representations from $PSL_2(\mathbb{Z})$ to
$GL_2(\textit{k})$ if $\textit{k}$ has characteristic 3. In fields
not of characteristic 3, we have
$$(a,b) = (1,\zeta), (1,\zeta^2), (\zeta,\zeta^2), (\zeta,1),
(\zeta^2, 1), \mbox{ or } (\zeta^2, \zeta).$$  Conjugating $X_1$ and
$Y_1$ by permutation matrix $X_1$ yields $(X_1,Y_1) \approx (X_1,
\left(
\begin{matrix} b & 0 \\ 0 & a \end{matrix} \right))$.  Thus we may
assume the following cases: $(a,b) = (1,\zeta), (1,\zeta^2)$, or
$(\zeta,\zeta^2)$. Since $tr(Y_1)$ is distinct for each of the
listed cases, each yields a separate equivalence class. This gives
the desired result.
\subsection{{n=3}{\bf}}
Let $\Lambda = \left(
\begin{matrix} \lambda_1 & & \\ & \lambda_2 &
\\ & & \lambda_3 \end{matrix} \right)$ and $\Gamma = \left(
\begin{matrix} \gamma_1 & & \\ & \gamma_2 & \\ & & \gamma_3
\end{matrix} \right)$. \smallskip By Lemma 3.4, we may assume that either $\lambda_1 =
\lambda_2 = \lambda_3$, $\lambda_1 = \lambda_2 \neq \lambda_3$, or
that $\lambda_1$, $\lambda_2$, and $\lambda_3$ are distinct. \\
\indent In the first case $\Lambda = \lambda_1 I$, which is a contradiction by Lemma 3.2. \\
\indent We will now solve the third case. By solving $\Lambda X
\Lambda = X$, we conclude for any $i, j \in \{1, 2, 3\}$ that
$\lambda_i \lambda_j = 1$ if $X_{i,j} \neq 0$. If $X$ has more than
one nonzero entry per row or column, we find that $\lambda_i =
\lambda_j$ for some $i \neq j$, which is a contradiction to our
assumption that $\lambda_1, \lambda_2, \lambda_3$ are distinct.
Therefore since $X$ is non-singular, $X$ has exactly one nonzero
entry per row and column. Solving $\Lambda Y \Gamma = Y$, we find
that for any $i, j \in \{1, 2, 3\}$, $\lambda_i \gamma_j = 1$ if
$Y_{i,j} \neq 0$. If $Y$ has more than one nonzero entry per row or
column, then either $\lambda_i = \lambda_j$ or $\gamma_i = \gamma_j$
for some $i \neq j$, both of which are contradictions. Therefore
since $Y$ is non-singular, $Y$ has exactly one nonzero entry per row
and column. Note, $\lambda_k = \frac{1}{\lambda_i}$ for some $i = 1,
2, 3$ and each $k$, by Remark 3.3. Since there are at most two
distinct solutions to $x = \frac{1}{x}$, then $\lambda_i \neq
\frac{1}{\lambda_i}$ for some $i$. Conjugating $X$ and $Y$ by a
permutation matrix if necessary, we may assume that $\lambda_2 =
\frac{1}{\lambda_1}$.  Thus $\lambda_3 = \frac{1}{\lambda_3}$.
Therefore, we find $\lambda_3 \lambda_1 \neq 1$ and $\lambda_3
\lambda_2 \neq 1$. Hence $X_{3,1} = X_{3,2} = 0$. This forces
$X_{3,3} \neq 0$. Expanding $\Gamma X \Gamma = X$, we see from
$X_{3,3} \neq 0$ that $\gamma_3 = \frac{1}{\gamma_3}$. Because
$\Lambda$ and $\Gamma$ are conjugate, $\gamma_3 = \lambda_3 =
\frac{1}{\lambda_3}$. Since $\gamma_3 = \frac{1}{\lambda_3}$ and
$\gamma_1, \gamma_2, \gamma_3$ are distinct, we find $\lambda_3
\gamma_2 \neq 1$ and $\lambda_3 \gamma_1 \neq 1$.  Thus $Y_{3,1} =
Y_{3,2} = 0$.  This forces $Y_{3,3} \neq 0$. But since $X$ and $Y$
have exactly one nonzero entry per row and column and $X_{3,3}$ and
$Y_{3,3}$ are nonzero, we find that $(X,Y)$ is reducible, which is a
contradiction.
\\
\indent In case 2, we assume $\lambda_1 = \lambda_2 \neq \lambda_3$.
First assuming $\lambda_1 \neq \frac{1}{\lambda_1}$ and using Remark
3.3, we see $\lambda_3 = \frac{1}{\lambda_1}$. Again solving
$\Lambda X \Lambda = X$, we conclude for any $i, j \in \{1, 2, 3\}$
that $\lambda_i \lambda_j = 1$ if $X_{i,j} \neq 0$. From $\lambda_3
= \frac{1}{\lambda_1}$ and $\lambda_1 = \lambda_2 \neq \lambda_3$,
we see that $\lambda_1^2 = \lambda_1 \lambda_2 = \lambda_2 \lambda_1
=\lambda_2^2 \neq \lambda_1 \lambda_3 = 1$. Then $X_{1,1} = X_{1,2}
= X_{2,1} = X_{2,2} = 0$. But this implies that $X$ is singular,
which is a contradiction. Hence $\lambda_1 = \frac{1}{\lambda_1}$
and $\lambda_3 = \frac{1}{\lambda_3}$.
\\
\indent In fields of characteristic 2, there is only one distinct
solution to $x=\frac{1}{x}$. This is a contradiction since $\lambda_1 \neq \lambda_3$.
Hence if $\textit{k}$ has characteristic
2, there are no irreducible representations from $PSL_2(\mathbb{Z})$
to $GL_3(\textit{k})$.
\\ \indent  We now consider
fields $\textit{k}$ which are not of characteristic 2. Since
$\lambda_i = \frac{1}{\lambda_i}$ for $i = 1, 2, 3$, then $\Lambda =
\pm \left(
\begin{matrix} 1 & & \\ & 1 & \\ & & -1 \end{matrix} \right)$.
Because $\Lambda$ and $\Gamma$ are conjugate and $\Lambda \neq
\Gamma$, we have 4 possibilities for $(\Lambda,\Gamma)$.
Substituting each of the possibilities for $(\Lambda, \Gamma)$ in
$\Lambda X \Lambda = X = \Gamma X \Gamma$, we solve the resulting
system of equations to find that in each possibility, $X$ is
diagonal. Again substituting each of the possibilities for
$(\Lambda, \Gamma)$ in the equation $\Lambda Y \Gamma = Y$ and using
$Y^3 = I$, we solve the resulting system of equations. In each
possibility, $Y$ has exactly one nonzero entry per row and
column.  \\
\indent Since $Y$ has exactly one nonzero entry per row and column,
\smallskip $Y^3 = I$, $X$ is
diagonal, and $(X,Y)$ is irreducible, we conclude $Y = \left(\begin{matrix} 0 & y_1 & 0 \\ 0 & 0 & y_2 \\
y_3 & 0 & 0 \end{matrix} \right)$ or $\left(\begin{matrix} 0 & 0 & y_1 \\ y_2 & 0 & 0 \\
0 & y_3 & 0 \end{matrix} \right)$, where $y_1y_2y_3 = 1$.  Conjugate $(X,Y)$ with $P = \left(
\begin{matrix} 1 & 0 & 0 \\ 0 & y1 & 0 \\ 0 & 0 & y_1y_2
\end{matrix} \right)$ in case 1, or with $P = \left(
\begin{matrix} 1 & 0 & 0 \\ 0 & 0 & y_1 \\ 0 & y_1y_3 & 0
\end{matrix} \right)$ in case 2, so that $(X,Y)
\approx (PXP^{-1},PYP^{-1})=:(X_1,Y_1)$, where $X_1$ is diagonal and
$Y_1 = \left(
\begin{matrix} 0 & 1 & 0 \\ 0 & 0 & 1 \\ 1 & 0 & 0 \end{matrix}
\right)$.  \smallskip Since $X_1^2 = I$, all of the diagonal entries
must be $\pm 1$.  Because $X_1 \neq \pm I$, there are 6
possibilities for $X_1$. Conjugating by various permutation
matrices, we can see that
$$(X,Y) \approx (\pm X_2, Y_1),\mbox{ where } X_2 = \left(
\begin{matrix} -1 & 0 & 0
\\ 0 & 1 & 0 \\ 0 & 0 & 1
\end{matrix} \right).$$
Thus $(X,Y)$ satisfies the hypotheses of Lemma 3.4 if and only if
$(\pm X_2, Y_1)$ satisfies the hypotheses of Lemma 3.4. Since
$(X_2,Y_1)$ satisfies $X_2^2 = I = Y_1^3$ and $X_2Y_1X_2Y_1^2$ and
$X_2Y_1^2X_2Y_1$ are diagonal, we need only check irreducibility. \\
\indent If $(X_2,Y_1)$ is irreducible, then $(-X_2,Y_1)$ is
irreducible. Hence it is sufficient to check irreducibility for
$(X_2,Y_1)$. However, we note that the standard basis matrices can
be composed as:
$$\frac{Y_1-Y_1X_2}{2}, \frac{Y_1(Y_1-Y_1X_2)}{2},
\frac{Y_1^2(Y_1-Y_1X_2)}{2},\frac{(Y_1-Y_1X_2)Y_1}{2},
\frac{Y_1(Y_1-Y_1X_2)Y_1}{2},$$
$$\frac{Y_1^2(Y_1-Y_1X_2)Y_1}{2}, \frac{(Y_1-Y_1X_2)Y_1^2}{2}, \frac{Y_1(Y_1-Y_1X_2)Y_1^2}{2},
 \mbox{ and } \frac{Y_1^2(Y_1-Y_1X_2)Y_1^2}{2}.$$ Note that $(X_2,Y_1)$ and $(-X_2,Y_1)$
yield separate equivalence classes since the trace of the matrices
are preserved in each equivalence class and $tr(X_2) \neq tr(-X_2)$.
We have therefore achieved the desired result.
\section{Case when n=6}
In this section, we prove (i), (ii), and (iii) of Theorem 4.1
in the case when $n=6$.
\begin{lemma}
Suppose $X, Y \in GL_n(\textit{k})$ where $X^2=Y^3=I$, $(X,Y)$ is
irreducible, and $\Lambda = XYXY^2$ and $\Gamma = XY^2XY$ are
diagonal matrices. Then $X$ and $Y$ must have exactly one nonzero
entry per row and column.
\end{lemma}
\begin{proof}
Assume $X$ and $Y$ satisfy the hypothesis of the lemma.  Let
$\Lambda$ and $\Gamma$ be diagonal matrices with diagonal entries
$\Lambda_{i,i} =: \lambda_i$, $\Gamma_{i,i} =: \gamma_i$. Note that
$(X,Y)$ satisfies the hypothesis if and only if
$(PXP^{-1},PYP^{-1})$ satisfies the hypothesis where $P$ is a
permutation matrix. Thus we may assume the diagonal entries of
$\Lambda$ are in one of the following cases, where entries in
distinct ordered tuples are not equal and entries in the same
ordered tuple are equal:
\begin{enumerate}
\begin{item}
$(\lambda_1, \lambda_2, \lambda_3, \lambda_4, \lambda_5, \lambda_6)$
\end{item}
\begin{item}
$(\lambda_1, \lambda_2, \lambda_3, \lambda_4, \lambda_5),
(\lambda_6)$
\end{item}
\begin{item}
$(\lambda_1, \lambda_2, \lambda_3, \lambda_4), (\lambda_5,
\lambda_6)$
\end{item}
\begin{item}
$(\lambda_1, \lambda_2, \lambda_3), (\lambda_4, \lambda_5,
\lambda_6)$
\end{item}
\begin{item}
$(\lambda_1, \lambda_2, \lambda_3, \lambda_4), (\lambda_5),
(\lambda_6)$
\end{item}
\begin{item}
$(\lambda_1, \lambda_2, \lambda_3), (\lambda_4, \lambda_5),
(\lambda_6)$
\end{item}
\begin{item}
$(\lambda_1, \lambda_2), (\lambda_3, \lambda_4), (\lambda_5,
\lambda_6)$
\end{item}
\begin{item}
$(\lambda_1, \lambda_2, \lambda_3), (\lambda_4), (\lambda_5),
(\lambda_6)$
\end{item}
\begin{item}
$(\lambda_1, \lambda_2), (\lambda_3, \lambda_4), (\lambda_5),
(\lambda_6)$
\end{item}
\begin{item}
$(\lambda_1, \lambda_2), (\lambda_3), (\lambda_4), (\lambda_5),
(\lambda_6)$
\end{item}
\begin{item}
$(\lambda_1), (\lambda_2), (\lambda_3), (\lambda_4), (\lambda_5),
(\lambda_6)$
\end{item}
\end{enumerate}
\underline{Case 1}: \\
In this case $\Lambda = \lambda_1 I$, which by Lemma 3.2 is a contradiction to the irreducibility of $(X,Y)$.\\
\\
\underline{Case 2, 3}: \\
First assume $\lambda_1 \neq \frac{1}{\lambda_1}$. Expanding
$\Lambda X \Lambda = X$ we get a system of equations where
for any $i, j \in \{1, \ldots, 6\}$, $\lambda_i\lambda_j = 1$ if $X_{i,j}\neq 0$. Since $\lambda_1 \neq
\frac{1}{\lambda_1}$, we find $X$ must be singular, which is a
contradiction. Thus assume $\lambda_1 = \frac{1}{\lambda_1}$. We
conclude from Remark 3.3 that $\lambda_6 = \frac{1}{\lambda_6}$ and
hence $\Lambda^2 = I$. Because $\Lambda \Gamma = \Gamma \Lambda$,
we see that $(XY)^6 = (YX)^6 = I$. Also, $(XY)^3 = (YX)^3$ since $\Lambda^2 =
I$. \\
\noindent Let $A=XY$ and $B=YX$. Since $A^3 = B^3$, $ABA =
BAB$, and $A^6 = B^6 = I$, we can count possible monomials in $A$
and $B$ to find that $(A,B)$ span at most a 24-dimensional space and
thus cannot be irreducible. Note $(A,B)$ is irreducible if and only
if $(X,Y)$ is irreducible, since we can generate $X$ and $Y$
from $A$ and $B$ and vice versa. Therefore $(X,Y)$ is not irreducible. \\
\\
\underline{Case 4}: \\
By Remark 3.3, if $\lambda_1 = \frac{1}{\lambda_1}$, then $\lambda_6
= \frac{1}{\lambda_6}$ and $\Lambda^2 = I$.  As in Cases 2 and 3,
this leads to a contradiction.  Hence $\lambda_6 =
\frac{1}{\lambda_1}$. Observe that since $\Lambda$ and $\Gamma$ are
conjugate,
\begin{equation}
\Lambda + \Lambda^{-1} = (\lambda_1 + \frac{1}{\lambda_1})I =\Gamma
+ \Gamma^{-1}.
\end{equation}
Then it follows that
\begin{equation}
Y(\Lambda + \Lambda^{-1}) = Y(\Gamma + \Gamma^{-1}) = (\Gamma +
\Gamma^{-1})Y.
\end{equation}
Substituting $\Lambda = XYXY^2$ and $\Gamma = XY^2XY$ into equation
(6.2) and using $(YX)^6=I$ yields
\begin{equation}
(YX)((YX)^2 - (XY)^2)(XY - YX) = 0.
\end{equation}
Since $XY = YX \Gamma$, we find that $(YX)((YX)^2 -
(XY)^2)(YX)(\Gamma-I) = 0.$ Given that $\Lambda$ and $\Gamma$ are
conjugate and $\lambda_1$ and $\lambda_6$ are not equal to 1,
$(\Gamma-I)$ is invertible. Also by invertibility of $(YX)$, we
deduce $(YX)^2 = (XY)^2$. Thus $A^2 = B^2$, $A^6 = B^6 = I$, and
$ABA = BAB$, where $A = XY$ and $B=YX$. Again by counting monomials
in $A$ and $B$, we find $A$ and $B$ span at most an 18-dimensional
space and thus $(A,B)$ cannot be irreducible. Since $(A,B)$ is
irreducible if and only if $(X,Y)$ is
irreducible, $(X,Y)$ cannot be irreducible. \\
\\
\underline{Case 5}: \\
First assume $\lambda_1 \neq \frac{1}{\lambda_1}$. Then by Remark
3.3, $\lambda_5$ or $\lambda_6$ must equal $\frac{1}{\lambda_1}$. We
expand $\Lambda X \Lambda = X$ to get a system of equations where
for any $i, j \in \{1, \ldots, 6\}$, $\lambda_i\lambda_j = 1$ if $X_{i,j}\neq 0$. If either $\lambda_5
= \frac{1}{\lambda_1}$ or $\lambda_6 = \frac{1}{\lambda_1}$, then
$X$ must be singular, which is a contradiction. Thus $\lambda_1 =
\frac{1}{\lambda_1}$.  If $\lambda_5 = \frac{1}{\lambda_5}$ then by
Remark 3.3, $\lambda_6 = \frac{1}{\lambda_6}$. This is a
contradiction since there are at most two distinct solutions to $x =
\frac{1}{x}$. Hence $\lambda_5 = \frac{1}{\lambda_6}$. Since
$\Lambda$ and
$\Gamma$ are conjugate, there are $ \left( \begin{array} {c} 6 \\
[1ex] 4,1,1 \end{array} \right)$ permutations of the $\lambda_i$'s.
Thus there are 30 possible matrices for $\Gamma$ with nonzero
entries determined by the $\lambda_i$'s. By case-by-case checking,
we see that the restrictions $\Gamma X \Gamma = X = \Lambda X
\Lambda$ and $X$ is non-singular leave 14 possible matrices for
$\Gamma$ (see Appendix 1). Note that $XY = YX$ if $\Lambda = \Gamma$
and $(XY)^2 = (YX)^2$ if $\Lambda \Gamma = I$. In both cases, we
find $(X,Y)$ is not irreducible as shown earlier. This leaves 12
possible
matrices for $\Gamma$.   \\
\indent Let $P$ be a permutation matrix. Note that $(X,Y)$ satisfies the
hypothesis of Lemma 6.1 if and only if $(PXP^{-1}, PYP^{-1})$ does.  Also, $X$ and $Y$
have one nonzero entry per row and column if and only if $PXP^{-1}$ and $PYP^{-1}$ do.
Thus for any permutation matrix $P$, we may replace $(X,Y)$ with $(PXP^{-1},PYP^{-1})$.
Further case-by-case checking gives us that for each of the
12 possible values of $\Gamma$, we are able to replace $(X,Y)$ with $(PXP^{-1},PYP^{-1})$
for an appropriate permutation matrix $P$, so that $\Lambda$ is
preserved, $\gamma_1 =$ $\gamma_2 =$ $\gamma_5 =$ $\gamma_6 =
\lambda_1$, $\gamma_3 = \lambda_5$, and $\gamma_4 = \lambda_6$ (see
Appendix 1). Substituting $\Lambda$ and $\Gamma$ into $\Lambda Y
\Gamma = Y$ and $\Gamma Y^2 \Lambda = Y^2$, using $Y^3 = I$, and solving
for $Y$, we find that $Y = \left(
\begin{matrix} 0 & 0 & Y_1 \\ Y_2 & 0 & 0
\\ 0 & Y_3 & 0 \end{matrix} \right)$, where $Y_1$, $Y_2$, and $Y_3$ are $2 \times 2$ block
matrices.  \smallskip Similarly, by substituting $\Lambda$ and
$\Gamma$ into $\Lambda X \Lambda = X = \Gamma X \Gamma$ and solving
for $X$, we find that \smallskip $X = \left(
\begin{matrix} X_1 & 0 & 0 \\ 0 & X_2 & 0 \\ 0 & 0 & X_3
\end{matrix} \right)$, where
$X_1$, $X_2$, and $X_3$ are $2 \times 2$ block matrices. \\
\indent Due to the fact that $Y^3 = I$, we find $Y_3 =
Y_1^{-1}Y_2^{-1}$. Computation with $XY = \Lambda YX$ yields
\begin{align*}
& X_1Y_1 = \lambda_1 Y_1 X_3 & & X_2Y_2 = \lambda_1 Y_2X_1 & X_3Y_3
= \left(
\begin{matrix} \lambda_5 & 0
\\ 0 & \frac{1}{\lambda_5} \end{matrix} \right) Y_3 X_2.
\end{align*}
But then:
$$X_3Y_3 = \frac{1}{\lambda_1} (Y_1^{-1}X_1Y_1)Y_3 =
\frac{1}{\lambda_1} (Y_1^{-1})(X_1Y_2^{-1}) = \frac{1}{\lambda_1}
(Y_1^{-1}) (\frac{1}{\lambda_1}) (Y_2^{-1} X_2) =
(\frac{1}{\lambda_1^2}) Y_3 X_2.$$  This is a contradiction to
$X_3Y_3 = \left(
\begin{matrix} \lambda_5 & 0
\\ 0 & \frac{1}{\lambda_5} \end{matrix} \right) Y_3 X_2$.
\\
\\
\underline{Case 6}: \\
Suppose $\lambda_6 \neq \frac{1}{\lambda_6}$.  Then either
$\lambda_1 = \frac{1}{\lambda_6}$ or $\lambda_4 =
\frac{1}{\lambda_6}$. We expand $\Lambda X \Lambda = X$ to get a
system of equations where for any $i, j \in \{1, \ldots, 6\}$, $\lambda_i\lambda_j = 1$ if $X_{i,j}\neq
0$. If either $\lambda_1 = \frac{1}{\lambda_6}$ or $\lambda_4 =
\frac{1}{\lambda_6}$, then $X$ must be singular, which is a
contradiction. Thus $\lambda_6 = \frac{1}{\lambda_6}$. Suppose next
that $\lambda_4 \neq \frac{1}{\lambda_4}$. Then $\lambda_1 =
\frac{1}{\lambda_4}$. Since $\lambda_i\lambda_j = 1$ if $X_{i,j}\neq
0$ and since $\lambda_1 = \frac{1}{\lambda_4}$, we find that $X$ is
singular, which is a contradiction.  Thus assume $\lambda_4 =
\frac{1}{\lambda_4}$. Because $\lambda_4 \neq \frac{1}{\lambda_1}$
and $\lambda_6 \neq \frac{1}{\lambda_1}$, it follows that $\lambda_1
= \frac{1}{\lambda_1}$. This is a contradiction since there are at
most 2 distinct solutions to $x = \frac{1}{x}$, while $\lambda_1$,
$\lambda_4$, and $\lambda_6$
are distinct. \\
\\
\underline{Case 7}: \\
It is easy to show using Remark 3.3 that there must be an
even number of ordered tuples with entries $\lambda_i$ such that
$\lambda_i \neq \frac{1}{\lambda_i}$. Since there are 3 tuples in
this case, there must be 1 or 3 tuples with entries $\lambda_i$
where $\lambda_i = \frac{1}{\lambda_i}.$ Because there are at most
two distinct solutions to $x = \frac{1}{x}$, there is only one tuple
with entries $\lambda_i$ where $\lambda_i = \frac{1}{\lambda_i}$.
Without loss of generality, assume that $\lambda_1 =
\frac{1}{\lambda_1}$.  Therefore $\lambda_3 = \frac{1}{\lambda_5}$.
Since $\Lambda$ and
$\Gamma$ are conjugate, there are $ \left( \begin{array} {c} 6 \\
[1ex] 2,2,2 \end{array} \right)= $ 90 possible matrices for $\Gamma$
in terms of the $\lambda_i$'s.  From $\Gamma X \Gamma = X = \Lambda
X \Lambda$, exactly 22 matrices for $\Gamma$ leave $X$ non-singular.
For each of these values of $\Gamma$, we find the corresponding form
for $X$ (see Appendix 1). Similarly, we substitute $\Lambda$ and
each possible $\Gamma$ into $\Lambda Y \Gamma = Y$ and solve the
resulting system of equations. For each of the 22 remaining values
of $\Gamma$, $Y$ must be of a certain corresponding form (See
Appendix 1). Further case-by-case analysis yields exactly 16 choices
of $\Gamma$ which do not force $(X,Y)$ to be reducible (see Appendix 1). In each of these 16
cases, the restriction that $XYXY^2$ is diagonal forces both $X$ and $Y$ to
have exactly one nonzero entry per row and column. \\
\\
\underline{Case 8}: \\
Suppose $\lambda_1 \neq \frac{1}{\lambda_1}$.  Then $\lambda_4$,
$\lambda_5$, or $\lambda_6$ must equal $\frac{1}{\lambda_1}$. We
expand $\Lambda X \Lambda = X$ to get a system of equations where
for any $i, j \in \{1, \ldots, 6\}$, $\lambda_i\lambda_j = 1$ if $X_{i,j}\neq 0$. If either $\lambda_1 =
\frac{1}{\lambda_4}$, $\lambda_1 = \frac{1}{\lambda_5}$, or
$\lambda_1 = \frac{1}{\lambda_6}$, then we find $X$ must be
singular, which is a contradiction.  Thus $\lambda_1 =
\frac{1}{\lambda_1}$. As in Case 7, it is easy to show that there
must be 0, 2, or 4 tuples with entries $\lambda_i$ where $\lambda_i
= \frac{1}{\lambda_i}$ for each $i$ in that tuple. Since $\lambda_1 =
\frac{1}{\lambda_1}$ and since there are at most 2 distinct
solutions to $x = \frac{1}{x}$, there must be two such tuples.
Without loss of generality, assume that $\lambda_4 =
\frac{1}{\lambda_4}$. (Note that if instead $\lambda_5 =
\frac{1}{\lambda_5}$ or $\lambda_6 = \frac{1}{\lambda_6}$, we may
conjugate $X$ and $Y$ by an appropriate permutation matrix P.) Since
$\lambda_5 \neq \frac{1}{\lambda_5}$, then $\lambda_6 =
\frac{1}{\lambda_5}$. Because $\Lambda$ and $\Gamma$ are conjugate,
there are $ \left( \begin{array} {c} 6 \\
[1ex] 3,1,1,1 \end{array} \right)= $ 120 possible matrices for
$\Gamma$ in terms of the $\lambda_i$'s. By case-by-case analysis, we
see from $\Gamma X \Gamma = X = \Lambda X \Lambda$ that exactly 20
choices of $\Gamma$ leave $X$ non-singular. For each of the 20
possible values of $\Gamma$, we find the corresponding form of $X$
(see Appendix 1). Then we substitute $\Lambda$ and each possible
$\Gamma$ into $\Lambda Y \Gamma = Y$ and solve the resulting system
of equations. For each possible $\Gamma$, $Y$ must be of a certain
corresponding form (see Appendix 1). Among these 20 choices of
$\Gamma$, further case-by-case analysis yields exactly 6 choices of
$\Gamma$ which do not force $(X,Y)$ to be reducible (see Appendix 1).  In each of the 6
remaining values of $\Gamma$, we find that $XYXY^2$ is
not diagonal, which is a contradiction. \\
\\
\underline{Case 9}: \\
Suppose that $\lambda_1 = \frac{1}{\lambda_1}$. Expand $\Lambda X
\Lambda = X$ to yield a system of equations where
for any $i, j \in \{1, \ldots, 6\}$, $\lambda_i\lambda_j = 1$ if $X_{i,j}\neq 0$. If $\lambda_3 \neq
\frac{1}{\lambda_3}$, then either $\lambda_5 = \frac{1}{\lambda_3}$
or $\lambda_6 = \frac{1}{\lambda_3}$. In either case, $X$ must be
singular, which is a contradiction. Thus $\lambda_3 =
\frac{1}{\lambda_3}$. Since there are at most two distinct solutions
to $x = \frac{1}{x}$, we find $\lambda_6 = \frac{1}{\lambda_5}$.
There are $ \left( \begin{array} {c} 6 \\
[1ex] 2,2,1,1 \end{array} \right)= $ 180 possible matrices for
$\Gamma$ defined in terms of the $\lambda_i$'s. Substitute $\Lambda$
and each possible $\Gamma$ into $\Gamma X \Gamma = X = \Lambda X
\Lambda$ and solve the resulting system of equations. By
case-by-case analysis, we see that exactly 20 choices of $\Gamma$
leave $X$ non-singular. For each of these 20 possible values of
$\Gamma$, we find the corresponding form of $X$ (see Appendix 1).
Similarly, substitute $\Lambda$ and each possible $\Gamma$ into
$\Lambda Y \Gamma = Y$ and solve the resulting system of equations.
For each of the 20 remaining possible values of $\Gamma$, we find
the corresponding form for $Y$ (see Appendix 1).  Among these values
of $\Gamma$, further case-by-case analysis yields exactly 4 choices
of $\Gamma$ which do not force $(X,Y)$ to be reducible (see Appendix 1). In each of these
cases, the corresponding forms of $X$ and $Y$ are $X = \left(
\begin{matrix} X_1 & 0 & 0 \\ 0 & X_2 & 0
\\ 0 & 0 & X_3 \end{matrix} \right)$ and $Y = \left(
\begin{matrix} 0 & 0 & Y_1 \\ Y_2 & 0 & 0 \\ 0 & Y_3 & 0
\end{matrix} \right)$ or $\left( \begin{matrix} 0 & Y_1 & 0 \\ 0 &
0 & Y_2 \\ Y_3 & 0 & 0 \end{matrix} \right)$, where $X_i$ and $Y_i$
are $2 \times 2$ block matrices.  Using a similar argument to that in Case 5, we
reach a contradiction with both forms of $(X,Y)$. \\
\indent Hence $\lambda_1 \neq \frac{1}{\lambda_1}$. We then see that
$\lambda_3$, $\lambda_5$, or $\lambda_6$ equals
$\frac{1}{\lambda_1}$.  Expand $\Lambda X \Lambda = X$ to get a
system of equations where for any $i, j \in \{1, \ldots, 6\}$, $\lambda_i\lambda_j = 1$ if $X_{i,j}\neq
0$. If either $\lambda_5 = \frac{1}{\lambda_1}$ or $\lambda_6 =
\frac{1}{\lambda_1}$, it follows that $X$ is singular, which is a
contradiction. Thus $\lambda_3 = \frac{1}{\lambda_1}$. Suppose
$\lambda_5 = \frac{1}{\lambda_5}$. Then $\lambda_6 =
\frac{1}{\lambda_6}$.  From $\Lambda X \Lambda = X$, one sees that
$X_{5,5}$ and $X_{6,6}$ are nonzero elements. Note $\gamma_5^2 =
\gamma_6^2 = 1$ since $\Gamma X \Gamma = X$. Because $\Lambda$ and
$\Gamma$ are conjugate, either $\gamma_5 = \lambda_5$ and $\gamma_6
= \lambda_6$, or $\gamma_5 = \lambda_6$ and $\gamma_6 = \lambda_5$.
In both cases, $(X,Y)$ is not irreducible by $\Lambda Y \Gamma = Y$,
which is a contradiction. Therefore $\lambda_3 =
\frac{1}{\lambda_1}$ and $\lambda_6 =
\frac{1}{\lambda_5}$.  There are $\left( \begin{array} {c} 6 \\
[1ex] 2,2,1,1 \end{array} \right) =$ 180 possible matrices for
$\Gamma$. By case-by-case analysis, we see from $\Gamma X \Gamma = X
= \Lambda X \Lambda$ that exactly 44 choices of $\Gamma$ leave $X$
non-singular. For each of the 44 possible values of $\Gamma$, we
find the corresponding form of $X$ (see Appendix 1). Substitute
$\Lambda$ and each of the possible values of $\Gamma$ into $\Lambda
Y \Gamma = Y$ and solve the resulting system of equations. For each
of the 44 possible values of $\Gamma$, we find the corresponding
form of $Y$ (see Appendix 1).  Among these possible values of
$\Gamma$, further case-by-case analysis yields exactly 32 choices of
$\Gamma$ which do not force $(X,Y)$ to be reducible (see Appendix 1).
For each of the 32 remaining choices of $\Gamma$, the fact that $XYXY^2$
is diagonal either eliminates the choice of $\Gamma$, or
forces both $X$ and $Y$ to have exactly one nonzero entry per row
and column. \\
\\
\underline{Case 10}: \\
Assume that $\lambda_1 \neq \frac{1}{\lambda_1}$.  Then either $\lambda_3$,
$\lambda_4$, $\lambda_5$, or $\lambda_6$ must equal
$\frac{1}{\lambda_1}$.  We expand $\Lambda X \Lambda = X$ to get a
system of equations where for any $i, j \in \{1, \ldots, 6\}$, $\lambda_i\lambda_j = 1$ if $X_{i,j}\neq
0$. If either $\lambda_1 = \frac{1}{\lambda_3}$,
$\lambda_1 = \frac{1}{\lambda_4}$, $\lambda_1 =
\frac{1}{\lambda_5}$, or $\lambda_1 = \frac{1}{\lambda_6}$, then
$X$ must be singular, which is a contradiction. Thus $\lambda_1 = \frac{1}{\lambda_1}$.
As in cases 7 and 8, it is easy to show that there must be 1, 3, or 5 tuples with entries $\lambda_i$ in which
$\lambda_i = \frac{1}{\lambda_i}$. Since there are at most 2
distinct solutions to $x = \frac{1}{x}$, there must be only one such
tuple.  Without loss of generality, we may assume that $\lambda_4 =
\frac{1}{\lambda_3}$ and thus $\lambda_6 = \frac{1}{\lambda_5}$.
Since $\Lambda$ and $\Gamma$ are conjugate, there are $ \left( \begin{array} {c} 6 \\
[1ex] 2,1,1,1,1 \end{array} \right)= $ 360 possible matrices for
$\Gamma$ in terms of the $\lambda_i$'s. By case-by-case analysis, we see that the restrictions $\Gamma X \Gamma = X
= \Lambda X \Lambda$ and $X$ is non-singular leave exactly 24 choices for $\Gamma$.
For each of the 24 possible values of $\Gamma$, we
find the corresponding form of $X$ (see Appendix 1). Substitute
$\Lambda$ and each of the possible values of $\Gamma$ into $\Lambda
Y \Gamma = Y$ and solve the resulting system of equations. For each
of the 24 possible values of $\Gamma$, we find the corresponding
form of $Y$ (see Appendix 1). Among these 24 choices of $\Gamma$, further
case-by-case analysis yields exactly 8 choices of $\Gamma$ which
do not force $(X,Y)$ to be reducible (see Appendix 1).  In each of
these 8 cases, the restriction that $XYXY^2$ is diagonal forces both
$X$ and $Y$ to have exactly one nonzero entry per row and column. \\
\\
\underline{Case 11}: \\
By solving $\Lambda X \Lambda = X$, we find that for any $i, j \in \{1, \ldots, 6\}$, $\lambda_i
\lambda_j = 1$ if $X_{i,j} \neq 0$. If $X$
has more than one nonzero entry per row or column, then
$\lambda_i = \lambda_j$ for some $i \neq j$, which is a
contradiction to our assumption that all $\lambda_i$ are distinct.
Therefore since $X$ is non-singular, $X$ has exactly one nonzero
entry per row and column. Solving $\Lambda Y \Gamma = Y$, we find
that for any $i, j \in \{1, \ldots, 6\}$, $\lambda_i \gamma_j = 1$ if $Y_{i,j} \neq 0$.
If $Y$ has more than one nonzero entry per row or
column, then either $\lambda_i = \lambda_j$ or $\gamma_i = \gamma_j$
for some $i \neq j$, both of which are contradictions. Therefore
since $Y$ is non-singular, $Y$ has exactly one nonzero entry per row
and column. \\
\end{proof}

\begin{lemma}
If $(X,Y)$ satisfies the hypotheses of Lemma 3.4 where $n=6$, then
$$(X,Y) \approx \Biggl( \left(
\begin{matrix} 0 & 1 & 0 & 0 & 0 & 0 \\ 1 & 0 & 0 & 0 & 0 & 0 \\ 0
& 0 & 0 & 1 & 0 & 0 \\ 0 & 0 & 1 & 0 & 0 & 0 \\ 0 & 0 & 0 & 0 & 0
& 1 \\ 0 & 0 & 0 & 0 & 1 & 0 \end{matrix} \right), \left(
\begin{matrix} 0 & 0 & 0 & 0 & c_1 & 0 \\ 0 & 0 & 0 & 0 & 0 & 1 \\
c_2 & 0 & 0 & 0 & 0 & 0 \\ 0 & 1 & 0 & 0 & 0 & 0 \\ 0 & 0 &
\frac{1}{c_1c_2} & 0 & 0 & 0 \\ 0 & 0 & 0 & 1 & 0 & 0 \end{matrix}
\right) \Biggr)$$ where $c_1,c_2 \in$ $\textit{k}$, $c_1,c_2 \neq
0$, and
\begin{enumerate}
\begin{item}
$(c_1,c_2) \neq (1,1), (-1,-1), (-1,1), (1,-1), (\zeta,\zeta),
(\zeta^2,\zeta^2)$, if $\textit{k}$ does not have characteristic 2
or 3, where $\zeta$ is a primitive cube root of unity.
\end{item}
\begin{item}
$(c_1,c_2) \neq (1,1), (\zeta,\zeta), (\zeta^2,\zeta^2)$, if
$\textit{k}$ has characteristic 2, where $\zeta$ is a primitive cube
root of unity.
\end{item}
\begin{item}
$(c_1,c_2) \neq (1,1), (-1,-1), (-1,1), (1,-1)$, if $\textit{k}$ has
characteristic 3.
\end{item}
\end{enumerate}
Furthermore, for $c_1, c_2$ satisfying the above criteria, the
following correspond to the same equivalence class: $(c_1,c_2)$,
$(\frac{1}{c_1}$,$\frac{1}{c_2})$, $(c_2,\frac{1}{c_1c_2})$,
$(\frac{1}{c_2},c_1c_2)$, $(\frac{1}{c_1c_2},c_1)$, and
$(c_1c_2,\frac{1}{c_1})$. If $(c'_1,c'_2)$ does not equal any
of the preceding tuples, then $(c'_1,c'_2)$ represents an equivalence class
distinct from $(c_1,c_2)$.
\end{lemma}
\begin{proof}
Assume $(X,Y)$ satisfies the hypotheses of Lemma 3.4. Let $H$ be the
group of $6 \times 6$ matrices with exactly one nonzero entry per
row and column and with the group operation of matrix
multiplication. Let $M$ be the group of all equivalence classes of
$H$ where $A \in H$ is equivalent to $B \in H$ when $$A_{i,j} = 0 \Leftrightarrow B_{i,j} = 0.$$  Denote $[A] \in M$ as the equivalence
class of $A \in H$ and denote $[1]_M$ as the identity element in
$M$. Now consider $[X] \in M$. Since $X^2 = I$, $[X]$ has order 1 or
2. We see that $M$ is isomorphic to $S_6$. The conjugacy classes of
$S_6$ with elements that have order 1 or 2 are represented by $
(1)$, $(1 \mbox{ } 2)$, $(1 \mbox{ } 2)$ $(3 \mbox{ } 4)$, and $(1
\mbox{ } 2)$$(3 \mbox{ } 4)$$(5 \mbox{ } 6)$. Thus there must be a
permutation matrix, $P$, such that $PXP^{-1}$ is diagonal or one of
\smallskip
\\ $X_0 = \left(
\begin{matrix} 0 & x_1& 0 & 0 & 0 & 0 \\ x_2 & 0 & 0 & 0 & 0 & 0
\\ 0 & 0 & x_3 & 0 & 0 & 0 \\ 0 & 0 & 0 & x_4 & 0 & 0 \\ 0 & 0 & 0 & 0 & x_5 & 0 \\ 0 & 0 & 0 & 0 & 0 & x_6
\end{matrix} \right)$, $X_1 = \left(
\begin{matrix} 0 & x_1& 0 & 0 & 0 & 0 \\ x_2 & 0 & 0 & 0 & 0 & 0
\\ 0 & 0 & 0 & x_3 & 0 & 0 \\ 0 & 0 & x_4 & 0 & 0 & 0 \\ 0 & 0 & 0 & 0 & x_5 & 0 \\ 0 & 0 & 0 & 0 & 0 & x_6
\end{matrix} \right)$, or \\ $X_2 = \left(
\begin{matrix} 0 & x_1& 0 & 0 & 0 & 0 \\ x_2 & 0 & 0 & 0 & 0 & 0
\\ 0 & 0 & 0 & x_3 & 0 & 0 \\ 0 & 0 & x_4 & 0 & 0 & 0 \\ 0 & 0 & 0 & 0 & 0 & x_5 \\ 0 & 0 & 0 & 0 &
x_6 & 0 \end{matrix} \right)$, for suitable choices of $x_1, \ldots, x_6$. \smallskip \\
Similarly, since $Y^3 = I$, there must be a permutation matrix, $Q$,
such that $QYQ^{-1}$ is diagonal or is equal to $$Y_1 = \left(
\begin{matrix} 0 & y_1& 0 & 0 & 0 & 0 \\ 0 & 0 & y_2 & 0 & 0 & 0
\\ y_3 & 0 & 0 & 0 & 0 & 0 \\ 0 & 0 & 0 & y_4 & 0 & 0 \\ 0 & 0 & 0 & 0 & y_5 & 0 \\ 0 & 0 & 0 & 0 & 0 & y_6
\end{matrix} \right) \mbox{ or } Y_2 = \left(
\begin{matrix} 0 & y_1& 0 & 0 & 0 & 0 \\ 0 & 0 & y_2 & 0 & 0 & 0
\\ y_3 & 0 & 0 & 0 & 0 & 0 \\ 0 & 0 & 0 & 0 & y_4 & 0 \\ 0 & 0 & 0 & 0 & 0 & y_5 \\ 0 & 0 & 0 & y_6 & 0 & 0
\end{matrix} \right)$$ for suitable choices of $y_1, \ldots, y_6$. \\
\indent Assume $PXP^{-1}$ is diagonal. Since $P \in H$, we find
$[1]_M = [PXP^{-1}] = [P][X][P^{-1}]$. Thus $[X] = [1]_M$, i.e. $X$
is a diagonal matrix.  Now take an arbitrary monomial in $X$ and
$Y$. Since $X^2 = I = Y^3$, this monomial can be expressed in the
form $Y^{a_1}XY^{a_2}X \ldots Y^{a_n}$ where $a_i = 0, 1, \mbox{ or
} 2$ for $i = 1 \ldots n$. Because $X$ is diagonal and $Y^3 = I$,
$$[Y^{a_1}XY^{a_2}X \ldots Y^{a_n}] = [Y^{a_1}][X][Y^{a_2}][X] \ldots
[Y^{a_n}] = [Y^{a_1}][1]_M[Y^{a_2}][1]_M \ldots [Y^{a_n}]$$ $$ =
[Y^{a_1}][Y^{a_2}] \ldots [Y^{a_n}] = [Y^{a_1+a_2+ \ldots + a_n}] =
[1]_M, [Y], \mbox{ or } [Y^2].$$  Since $Y \in H$ and $Y^2 \in H$,
we find that $Y$ and $Y^2$ have exactly one nonzero entry per row
and column. Thus $(X,Y)$ can span at most an 18-dimensional
space, which is a contradiction to $(X,Y)$ being irreducible. \\
\indent Now assume $QYQ^{-1}$ is diagonal for some permutation
matrix $Q$. By a similar argument to that above, we see that $Y$ is
a diagonal matrix. Now take an arbitrary monomial in $X$ and $Y$.
Again, this monomial can be expressed as $Y^{a_1}XY^{a_2}X \ldots Y^{a_n}$
where $a_i = 0, 1, \mbox{ or } 2$ for $i = 1 \ldots n$. Then since
$Y$ is diagonal and $X^2 = I$,
$$[Y^{a_1}XY^{a_2}X \ldots Y^{a_n}] = [Y^{a_1}][X][Y^{a_2}][X] \ldots
[Y^{a_n}] = [1]_M[X][1]_M[X] \ldots [1]_M$$ $$ = [X][X] \ldots [X] =
[X^{n-1}] = [1]_M \mbox{ or } [X].$$ Since $X \in H$, we find
$X$ has one nonzero entry per row and column. Therefore $(X,Y)$ can
span at most a 12-dimensional space, which is a contradiction to
$(X,Y)$ being irreducible. Hence $PXP^{-1} = X_0, X_1, \mbox{ or }
X_2$ and $QYQ^{-1} = Y_1$ or $Y_2$, for some permutation matrix $Q$.  We define $Y_0$ to be
$PYP^{-1}$. Then for each $(X,Y)$ which satisfies the hypotheses of Lemma 3.4,
we find that $(X,Y)$ must be equivalent to $(X_0,Y_0)$, $(X_1,Y_0)$, or $(X_2,Y_0)$ where
$Y_0$ is conjugate to $Y_1$ or $Y_2$ by a permutation matrix. Note
that if $P$ is a weighted permutation matrix, then $(X,Y)$ satisfies
the hypotheses of Lemma 3.4 if and only if $(PXP^{-1},PYP^{-1})$
satisfies the hypotheses of Lemma 3.4. Thus $(X_i,Y_0)$
satisfies the hypotheses of Lemma 3.4.  \\
\indent After examining all possible pairs $(X_i,Y_0)$ by using the
restrictions that $(X_i,Y_0)$ is irreducible and $X_iY_0X_iY_0^2$ is diagonal, we find
$(X,Y) \approx (X_2,Y_0)$, where $X_2$ is of the form above and $Y_0$ equals:
$$\left(
\begin{matrix} 0 & 0 & y_1 & 0 & 0 & 0 \\ 0 & 0 & 0 & y_2 & 0 &
0 \\ 0 & 0 & 0 & 0 & y_3 & 0 \\ 0 & 0 & 0 & 0 & 0 & y_4 \\ y_5 & 0 &
0 & 0 & 0 & 0 \\ 0 & y_6 & 0 & 0 & 0 & 0 \end{matrix} \right),
\left( \begin{matrix} 0 & 0 & y_1 &
0 & 0 & 0 \\ 0 & 0 & 0 & y_2 & 0 & 0 \\ 0 & 0 & 0 & 0 & 0 & y_3 \\
0 & 0 & 0 & 0 & y_4 & 0 \\ 0 & y_5 & 0 & 0 & 0 & 0 \\ y_6 & 0 &
0 & 0 & 0 & 0 \end{matrix} \right), \left( \begin{matrix} 0 & 0 & 0&
y_1 & 0 & 0 \\ 0 & 0 & y_2 & 0 & 0 & 0 \\ 0 & 0 & 0 & 0 & y_3 & 0
\\ 0 & 0 & 0 & 0 & 0 & y_4 \\ 0 & y_5 & 0 & 0 & 0 & 0 \\ y_6 & 0 &
0 & 0 & 0 & 0 \end{matrix} \right),$$ $$\left( \begin{matrix} 0 &
0 & 0 & y_1 & 0 & 0 \\ 0 & 0 & y_2 & 0 & 0 & 0 \\ 0 & 0 & 0 & 0 &
0 & y_3 \\ 0 & 0 & 0 & 0 & y_4 & 0 \\ y_5 & 0 & 0 & 0 & 0 & 0 \\ 0
& y_6 & 0 & 0 & 0 & 0 \end{matrix} \right), \left( \begin{matrix}
0 & 0 & 0 & 0 &
0 & y_1 \\ 0 & 0 & 0 & 0 & y_2 & 0 \\ y_3 & 0 & 0 & 0 & 0 & 0 \\
0 & y_4 & 0 & 0 & 0 & 0 \\ 0 & 0 & 0 & y_5 & 0 & 0 \\ 0 & 0 & y_6
& 0 & 0 & 0 \end{matrix} \right), \left( \begin{matrix} 0 & 0 & 0
& 0 & y_1 & 0 \\ 0 & 0 & 0 & 0 & 0 & y_2 \\ y_3 & 0 & 0 & 0 & 0 &
0 \\ 0 & y_4 & 0 & 0 & 0 & 0 \\ 0 & 0 & y_5 & 0 & 0 & 0 \\ 0 & 0 &
0 & y_6 & 0 & 0 \end{matrix} \right),$$ $$\left( \begin{matrix} 0
& 0 & 0 & 0 & y_1 & 0 \\ 0 & 0 & 0 & 0 & 0 & y_2 \\ 0 & y_3 & 0 &
0 & 0 & 0 \\ y_4 & 0 & 0 & 0 & 0 & 0 \\ 0 & 0 & 0 & y_5 & 0 & 0 \\
0 & 0 & y_6 & 0 & 0 & 0 \end{matrix} \right),  \mbox{ or } \left(
\begin{matrix} 0 & 0 & 0 & 0 & 0 & y_1 \\ 0 & 0 & 0 & 0 & y_2 & 0
\\ 0 & y_3 & 0 & 0 & 0 & 0 \\ y_4 & 0 & 0 & 0 & 0 & 0 \\ 0 & 0 &
y_5 & 0 & 0 & 0 \\ 0 & 0 & 0 & y_6 & 0 & 0 \end{matrix} \right),$$
for suitable choices of $y_1, \ldots, y_6$. From computation with
$X^2 = Y^3 = I$ and \smallskip conjugation by various weighted
permutation matrices, one sees that $$(X,Y) \approx (X',Y') :=
\Biggl( \left(
\begin{matrix} 0 & 1 & 0 & 0 & 0 & 0 \\ 1 & 0 & 0 & 0 & 0 & 0 \\ 0
& 0 & 0 & 1 & 0 & 0 \\ 0 & 0 & 1 & 0 & 0 & 0 \\ 0 & 0 & 0 & 0 & 0
& 1 \\ 0 & 0 & 0 & 0 & 1 & 0 \end{matrix} \right), \left(
\begin{matrix} 0 & 0 & 0 & 0 & c_1 & 0 \\ 0 & 0 & 0 & 0 & 0 & 1 \\
c_2 & 0 & 0 & 0 & 0 & 0 \\ 0 & 1 & 0 & 0 & 0 & 0 \\ 0 & 0 &
\frac{1}{c_1c_2} & 0 & 0 & 0 \\ 0 & 0 & 0 & 1 & 0 & 0 \end{matrix}
\right) \Biggr )$$ for some nonzero $c_1, c_2 \in$ $\textit{k}$.
\\
\indent We note $X'^2 = I = Y'^3$, $X'Y'X'Y'^2$ is diagonal, and
$X'Y'^2X'Y'$ is diagonal. Thus it remains to determine which values
of $c_1$ and $c_2$ lead to an irreducible $(X',Y')$ and to find the
equivalence classes for $(X',Y')$. We first determine the
equivalence classes for solutions
$(X',Y').$ \smallskip \\
Let $Y^* = \left( \begin{matrix} 0 & 0 & 0 & 0 & c'_1 & 0 \\
0 & 0 & 0 & 0 & 0 & 1 \\ c'_2 & 0 & 0 & 0 & 0 & 0 \\ 0 & 1 & 0 & 0 &
0 & 0
\\ 0 & 0 & \frac{1}{c'_1c'_2} & 0 & 0 & 0 \\ 0 & 0 & 0 & 1 & 0 & 0
\end{matrix} \right)$. \smallskip  Then $(X',Y^*) \approx (X',Y')$ if and only if there
exists an invertible matrix $Q$ such that $QX' = X'Q$ and
$QY'=Y^*Q.$ \smallskip \\ Let $Q = \left( \begin{matrix} A & B & C \\ D & E & F
\\ G & H & J \end{matrix} \right)$ where $A$ through $J$ are $2
\times 2$ block matrices. Then by $QX' = X'Q$, we find that $A = \left(
\begin{matrix} a_1 & a_2 \\ a_2 & a_1 \end{matrix} \right)$, $B = \left(
\begin{matrix} b_1 & b_2 \\ b_2 & b_1 \end{matrix} \right)$,
etc. \smallskip\\  From $QY' = Y^*Q$, we see that:
\begin{align*}
& A=0
\Leftrightarrow E=0 \Leftrightarrow J=0,& & B=0 \Leftrightarrow F=0
\Leftrightarrow G=0, & C=0 \Leftrightarrow D=0
\Leftrightarrow H=0.
\end{align*}
From this result and again using $QY' =
Y^*Q$, we find $(c'_1,c'_2)=(c_1,c_2) \mbox{ or }
(\frac{1}{c_1},\frac{1}{c_2})$ if $A \neq 0$. Likewise,
$(c'_1,c'_2)=(c_2,\frac{1}{c_1c_2}) \mbox{ or }
(\frac{1}{c_2},c_1c_2)$ if $B \neq 0$, and
$(c'_1,c'_2)=(\frac{1}{c_1c_2},c_1) \mbox{ or }
(c_1c_2,\frac{1}{c_1})$ if $C \neq 0$.\smallskip \\
Note that using:
\begin{align*}
& Q_1 = I,& & Q_2 = \left(
\begin{matrix} 0 & 1 &  &  &  &  \\ 1 & 0 &  &  &  &  \\ &  & 0 &
\frac{1}{c_2} &  &  \\ & & \frac{1}{c_2} & 0 & & \\ & & & & 0 & c_1 \\
& & & & c_1 & 0
\end{matrix} \right),  \\ & Q_3 = \left( \begin{matrix} & & 1 & 0 &
& \\ & & 0 & 1 & & \\ & & & & 1 & 0 \\ & & & & 0 & 1 \\ 1 & 0 & & &
& \\ 0 & 1 & & & & \end{matrix} \right), & & Q_4 =
\left(
\begin{matrix} & & 0 & 1 & & \\ & & 1 & 0 & & \\ & & & & 0 & c_1c_2 \\
& & & & c_1c_2 & 0 \\ 0 & c_2 & & & & \\ c_2 & 0 & & & &
\end{matrix} \right), \\ & Q_5 = \left( \begin{matrix} & & & & 1 & 0
\\ & & & & 0 & 1 \\ 1 & 0 & & & & \\ 0 & 1 & & & & \\ & & 1 & 0 &
& \\ & & 0 & 1 & & \end{matrix} \right), & & Q_6 = \left(
\begin{matrix} & & & & 0 & c_1c_2 \\ & & & & c_1c_2 & 0 \\ 0 & c_2
& & & & \\ c_2 & 0 & & & & \\ & & 0 & 1 & & \\ & & 1 & 0 & &
\end{matrix} \right)
\end{align*}
we find $Q_i X' Q_i^{-1} = X'$ and $Q_i Y' Q_i^{-1} = Y^*$ for $i = 1 \ldots 6$,
and $(c'_1,c'_2)$ as above.  Since $Q$ is invertible, $A$, $B$, $C$
cannot all be 0. Thus $Y^*$ must have
$$(c'_1,c'_2) = (c_1,c_2), (\frac{1}{c_1}, \frac{1}{c_2}), (c_2,
\frac{1}{c_1c_2}), (\frac{1}{c_2}, c_1c_2),
(\frac{1}{c_1c_2}, c_1), \mbox{ or } (c_1c_2, \frac{1}{c_1}).$$ \\
It remains to check irreducibility of $(X',Y').$ \smallskip \\
\indent Note that conjugating by $P = \left( \begin{matrix} 1 & 0 & 0 & 0 & 0 & 0 \\
0 & 0 & 0 & 0 & 0 & 1 \\ 0 & 0 & 1 & 0 & 0 & 0 \\ 0 & 1 & 0 & 0 &
0 & 0 \\ 0 & 0 & 0 & 0 & 1 & 0 \\ 0 & 0 & 0 & 1 & 0 & 0
\end{matrix} \right)$ yields $(X',Y') \approx (X'',Y')$ \smallskip where $X'' = \left(
\begin{matrix} 0 & 0 & 0 & 1 & 0 & 0 \\ 0 & 0 & 0 & 0 & 1 & 0 \\ 0
& 0 & 0 & 0 & 0 & 1
\\ 1 & 0 & 0 & 0 & 0 & 0 \\ 0 & 1 & 0 & 0 & 0 & 0 \\ 0 & 0 & 1 & 0
& 0 & 0 \end{matrix} \right)$. Thus $(X',Y')$ is irreducible if and
only if $(X'',Y')$ is irreducible. \smallskip \\
\indent Let $L = \left(
\begin{matrix} 0 & 0 & 0 & 0 & 0 & 0 \\ 0 & 0 & 0 & 0 & 0 & 0 \\ 0
& 0 & 0 & 0 & 0 & 0 \\ 0 & 0 & 0 & 0 & 0 & 0 \\ 0 & 0 & 0 & 0 & 0 &
0 \\ 1 & 0 & 0 & 0 & 0 & 0 \end{matrix} \right)$ and $U = \left(
\begin{matrix} 0 & 1 & 0 & 0 & 0 & 0 \\ 0 & 0 & 1 & 0 & 0 & 0 \\ 0
& 0 & 0 & 1 & 0 & 0 \\ 0 & 0 & 0 & 0 & 1 & 0 \\ 0 & 0 & 0 & 0 & 0 &
1 \\ 0 & 0 & 0 & 0 & 0 & 0 \end{matrix} \right)$.  \smallskip Since
$(L,U)$ is irreducible, it is sufficient to generate $L$ and $U$ as
$\textit{k}$-linear combinations of products in $X''$,
$Y'$.  Note that this is a necessary and sufficient condition for $(X'',Y')$ to be irreducible. \\
\indent We now find conditions that let us generate $L$ and $U$ from
$X''$ and $Y'$. Let $\Lambda'' = X''Y'X''Y'^2$ and $\Gamma'' =
X''Y'^2X''Y'$.
\smallskip Since $Y'X''$, $\Lambda'' Y'X''$, $\Lambda'' \Gamma'' Y'X''$,
$\Gamma'' Y'X''$, $\Lambda''^2 Y'X''$, and $\Gamma''^2 Y'X''$ are
matrices of the form \smallskip $\left(
\begin{matrix} 0 & * & 0 & 0 & 0 & 0 \\ 0 & 0 & * & 0 & 0 & 0 \\ 0
& 0 & 0 & * & 0 & 0 \\ 0 & 0 & 0 & 0 & * & 0 \\ 0 & 0 & 0 & 0 & 0 &
* \\ * & 0 & 0 & 0 & 0 & 0 \end{matrix} \right)$, we attempt to form
$U$ and $L$ as a $\textit{k}$-linear combination of those matrices.
For $a_i \in$ $\textit{k}$, \\ $$a_1Y'X'' + a_2 \Lambda'' Y'X'' +
a_3 \Lambda'' \Gamma'' Y'X''+a_4 \Gamma'' Y'X'' + a_5 \Lambda''^2
Y'X'' + a_6 \Gamma''^2 Y'X'' = U$$
\smallskip if and only if $A_1 \left(
\begin{matrix} a_1 \\ a_2 \\ a_3 \\ a_4 \\ a_5 \\ a_6 \end{matrix} \right)
=\left(
\begin{matrix} 0 \\ 1 \\ 1 \\ 1 \\ 1 \\ 1 \end{matrix} \right)$
\smallskip where $$A_1 = \left( \begin{matrix} 1 & c_2 & c_1c_2^2 &
c_1c_2 & c_2^2 & c_1^2c_2^2 \\ c_1 & 1 & c_2 & c_1c_2 &
\frac{1}{c_1} & c_1c_2^2 \\ 1 & \frac{1}{c_1c_2} &
\frac{1}{c_1^2c_2} & \frac{1}{c_1} & \frac{1}{c_1^2c_2^2} &
\frac{1}{c_1^2} \\ c_2 & 1 & \frac{1}{c_1c_2} & \frac{1}{c_1} &
\frac{1}{c_2} & \frac{1}{c_1^2c_2} \\ 1 & c_1 & \frac{c_1}{c_2} &
\frac{1}{c_2} & c_1^2 & \frac{1}{c_2^2} \\ \frac{1}{c_1c_2} & 1 &
c_1 & \frac{1}{c_2} & c_1c_2 & \frac{c_1}{c_2} \end{matrix}
\right).$$ \\  Similarly, for $b_i \in$ $\textit{k}$, \\ $$b_1Y'X''
+ b_2 \Lambda'' Y'X'' + b_3 \Lambda'' \Gamma'' Y'X''+b_4 \Gamma''
Y'X'' + b_5 \Lambda''^2 Y'X'' + b_6 \Gamma''^2 Y'X'' = L$$
\smallskip if and only if $A_1 \left(
\begin{matrix} b_1 \\ b_2 \\ b_3 \\ b_4 \\ b_5 \\ b_6 \end{matrix} \right)
=\left(
\begin{matrix} 1 \\ 0 \\ 0 \\ 0 \\ 0 \\ 0 \end{matrix} \right)$.
Thus $U$ and $L$ can be generated as above if
$\det(A_1)=R_1/(c_1^6c_2^6) \neq 0$, where $R_1$ is a polynomial in
$c_1$ and $c_2$. Likewise, we attempt to generate $U$ and $L$ using
matrices $Y'X''$, $Y'X'' \Lambda''$, $\Lambda'' Y'X''$, $\Gamma''
Y'X''$, $\Lambda'' \Gamma'' Y'X''$, and $\Lambda''^2Y'X''$.  Using
the method above, we obtain a matrix
$$B_1 = \left( \begin{matrix} 1 & \frac{1}{c_1} & c_2 & c_1c_2 &
c_1c_2^2 & c_2^2
\\ c_1 & \frac{1}{c_2} & 1 & c_1c_2 & c_2 & \frac{1}{c_1} \\ 1 &
\frac{1}{c_2} & \frac{1}{c_1c_2} & \frac{1}{c_1} &
\frac{1}{c_1^2c_2} & \frac{1}{c_1^2c_2^2} \\ c_2 & c_1c_2 & 1 &
\frac{1}{c_1} & \frac{1}{c_1c_2} & \frac{1}{c_2}
\\ 1 & c_1c_2 &
c_1 & \frac{1}{c_2} & \frac{c_1}{c_2} & c_1^2 \\
\frac{1}{c_1c_2} & \frac{1}{c_1} & 1 & \frac{1}{c_2} & c_1 & c_1c_2
\end{matrix} \right).$$  We are able to generate $L$ and $U$ from
the monomials above if $\det(B_1)=R_2/(c_1^5c_2^5) \neq 0$, where
$R_2$ is a polynomial in $c_1$ and $c_2$.  Thus, if we are not able
to generate
$L$ and $U$ in either way, $R_1 = R_2 = 0$. \\
\indent Earlier, we concluded that $(X',Y')$ is irreducible if and
only if $(X'',Y')$ is irreducible. Also, $(X'',Y')$ is irreducible
if and only if $L$ and $U$ can be generated from $X''$ and $Y'$.
From this, we are able to conclude that if
$(X',Y')$ is not irreducible, then $R_1=R_2=0$. \\
\indent We now solve $R_1 = R_2 = 0$. By factoring $R_2$ we find
$$R_2 = (1-c_1)(c_2-1)(c_1c_2-1)F_1F_2F_3$$ where $F_1$, $F_2$,
$F_3$ are irreducible polynomials in $c_1$, $c_2$. Our explicit choices of $F_1$, $F_2$, and
$F_3$ can be found in Appendix 2. Then $R_1 = R_2 = 0$ if and
only if $R_1 = 0$ and at least one of the following conditions hold:
$c_1 = 1$, $c_2 = 1$, $c_1c_2 = 1$, $F_1 = 0$, $F_2 = 0$, or $F_3 =
0$.
\begin{itemize}
\begin{item}
If $R_1 = 0$ and $c_1 = 1$ then $(c_2-1)^6(c_2^2-1)^3 = 0$.
\end{item}
\begin{item}
If $R_1 = 0$ and $c_2 = 1$ then $(c_1-1)^6(c_1^2-1)^3 = 0$.
\end{item}
\begin{item}
If $R_1 = 0$ and $c_1c_2 = 1$ then $(c_2-1)^6(c_2^2-1)^3 = 0$.
\end{item}
\end{itemize}
 Thus if $R_1=0$ and either $c_1=1$, $c_2 = 1$, or $c_1c_2 = 1$, then
$(c_1,c_2) = (1,1)$, $(1,-1)$, $(-1,1)$, or $(-1,-1)$ in fields not
of characteristic 2 and $(c_1,c_2) = (1,1)$ in fields of
characteristic 2. \\
\indent We next repeatedly use the command sprem in Maple to find a
univariate polynomial in $c_1$ or $c_2$ which is contained in the
ideal $\langle R_1,F_i \rangle$ of $\textit{k}[c_1,c_2]$ for each $i
= 1, 2, 3$. See Appendix 2 for the detailed Maple commands used in
this section of the paper. Specifically, sprem inputs a variable
$x$, and multivariate polynomials in $x$, say $a$ and $b$. It
computes multivariate polynomials in $x$ with integer coefficients,
say $m$ and $q$, where $m a = b q + r$, and the degree of $x$ in $r$
is strictly less than the degree of $x$ in $b$.  The output of sprem
is the multivariate function $r$.  $m$ is always of the form $x^n$
for some $n$. Temporarily regard $c_1$ and $c_2$ as indeterminates,
and temporarily replace $R_1$ and $F_1$ with their natural preimages
in $\mathbb{Z}[c_1,c_2]$. Using the function sprem to recursively
reduce the degree of $c_2$ when starting with initial polynomials
$R_1$ and $F_1$, we see that $c_1^{34} (c_1-1)^{24} (c^2_1+c_1+1)^6$
is in the ideal $\langle R_1,F_1 \rangle$ of $\mathbb{Z}[c_1,c_2]$.
Now let $c_1$, $c_2$, $R_1$, and $F_1$ again be elements of $\textit{k}$.
It follows that $c_1^{34} (c_1-1)^{24} (c^2_1+c_1+1)^6$ is in the ideal $\langle R_1,F_1 \rangle$ of
$\textit{k}[c_1,c_2]$. Hence if $R_1 = F_1 = 0$, either
$c_1 = 0$, $c_1 = 1$, or $c^2_1+c_1+1=0$. \\
\indent Similarly, we temporarily regard $c_1$ and $c_2$ as
indeterminates, and temporarily replace $R_1$ and $F_2$ with their natural preimages in
$\mathbb{Z}[c_1,c_2]$. We use the Maple command sprem to recursively reduce the
degree of $c_2$ when starting with initial polynomials $R_1$ and
$F_2$.  Our computation yield that $(c_1-1)^{24} (c^2_1+c_1+1)^6
(c^2_1 - c_1 + 1)^{11}$ is in the ideal $\langle R_1,F_2 \rangle$ of $
\mathbb{Z}[c_1,c_2]$. Returning $c_1, c_2, R_1,$ and $F_2$ to their original form, we conclude that $(c_1-1)^{24} (c^2_1+c_1+1)^6
(c^2_1 - c_1 + 1)^{11}$ is in the ideal $\langle R_1,F_2 \rangle$ of $\textit{k}[c_1,c_2]$. Hence if $R_1 = F_2 = 0$, either $c_1
= 1$, $c^2_1 + c_1 +1 = 0$, or $c^2_1 - c_1 + 1 = 0$. \\
\indent Again, temporarily regard $c_1$ and $c_2$ as
indeterminates, and temporarily replace $R_1$ and $F_3$ with their natural preimages in $\mathbb{Z}[c_1,c_2]$. We now solve $R_1 = F_3 = 0$ in two different ways.  First,
we recursively reduce the degree of $c_2$ by using the Maple command
sprem with initial polynomials $R_1$ and $F_3$.  We find
$$c_1^{186} (c_1+1)^{36} (c_1-1)^{56} (c_1^2+c_1+1)^{57} T_1^4
T_2^2$$ is in the ideal $\langle R_1,F_3 \rangle \subset \mathbb{Z}[c_1,c_2]$
where $T_1$ is an irreducible polynomial in $c_1$ of degree 28 and
where $T_2$ is an irreducible polynomial in $c_1$ of degree 40. Alternatively,
factor $R_1$ into $(c_2-c_1)(c_1^2c_2 - 1)(c_1c_2^2 - 1) R_3$,
where $R_3$ is a nonhomogeneous polynomial in $c_1$ and $c_2$ of
total degree 14.  Now we use the Maple command sprem to recursively reduce
the degree of $c_2$ using initial polynomials $F_3$ and each of the
polynomials $(c_2-c_1)$, $(c_1^2c_2 - 1)$, $(c_1c_2^2 - 1)$, and
$R_3$.  In respective order, we find that the following polynomials
are in $\langle R_1,F_3 \rangle \subset \mathbb{Z}[c_1,c_2]$:
\begin{itemize}
\begin{item}
$(c_1+1)^2 (c_1-1)^2 (c_1^2+c_1+1)^2$
\end{item} \\
\begin{item}
$c_1^2 (c_1+1)^2 (c_1-1)^2 (c_1^2+c_1+1)^2$
\end{item} \\
\begin{item}
$c_1^5 (c_1-1)^4 (c_1^2+c_1+1)^2$
\end{item} \\
\begin{item}
$c_1^{110} (c_1+1)^{16} (c_1-1)^{44} (c_1^2+c_1+1)^{31} T_3^4 T_4^2$
\end{item}
\end{itemize}
where $T_3$ is an irreducible polynomial of degree 16 and $T_4$ is
an irreducible polynomial of degree 26. Now remove the temporary replacements of $c_1$,
$c_2$, $R_1$, and $F_3$. By comparing solutions obtained in the two different ways,
we conclude that if $R_1 = F_3 = 0$, then $c_1 = 0$, $c_1-1 = 0$, $c_1 + 1 = 0$, or $c_1^2 + c_1 + 1 = 0$. \\
\indent Thus if $R_1 = 0$ and either $F_1$, $F_2$, or $F_3 = 0$,
then $c_1 = 0$, $c_1 = 1$, $c_1+1=0$, $c_1^2+c_1+1=0$, or
$c_1^2-c_1+1=0$. We reject the case that $c_1 = 0$ since this forces
$Y^3 \neq I$. The case where $c_1 = 1$ was solved above. If
$\textit{k}$ is not of characteristic 2, solving $c_1+1=0$ gives
$c_1 = -1$. If $c_1 = -1$ and $R_1 = R_2 = 0$, then $c_2 = 1$ or
$c_2 = -1$. These solutions are both listed above. If $\textit{k}$ is of
characteristic 2 and $c_1+1=0$ or $\textit{k}$ is of characteristic
3 and $c_1^2+c_1+1=0$ then $c_1 = 1$, which was solved above.  Also,
if $\textit{k}$ is of characteristic 3 and $c_1^2-c_1+1=0$ then $c_1
= -1$, which was solved above. It remains to consider fields not of
characteristic 3 where $c_1^2+c_1+1=0$ or $c_1^2-c_1+1=0$.
\\
\indent Thus assume $\textit{k}$ is not of characteristic 3,
$c_1^2+c_1+1=0$, and $R_1 = R_2 = 0$. Temporarily replace $c_1$, $c_2$, $R_1$, and $R_2$
with their natural preimages in $\mathbb{Z}[c_1,c_2]$. We use the Maple command sprem to recursively
reduce the degree of $c_1$ with initial polynomials $R_1$
and $c_1^2+c_1+1$, and $R_2$ and $c_1^2+c_1+1$.
We find that: \\
$$(c_2^2-c_2+1) (c_2^6-5c_2^5+23c_2^4-8c_2^3-c_2^2-2c_2+1)(c_2^2+c_2+1)^5 $$
$$(c_2^6-2c_2^5-c_2^4-8c_2^3+23c_2^2-5c_2+1) \in
\langle R_1, c_1^2+c_1+1 \rangle, \mbox{ and} \smallskip$$
$$27c_2^2(c_2^2-2c_2+4)(4c_2^2-2c_2+1)(c_2-1)^2(c_2^2+c_2+1)^5 \in \langle R_2,
c_1^2+c_1+1 \rangle.$$ \\ Note that if $R_1= R_2= 0$ and $c_1^2+c_1+1=0$,
then $c_2^2 + c_2 + 1 = 0$. Now remove the temporary replacements of $c_1$,
$c_2$, $R_1$, and $R_2$. Thus $c_1 = \zeta \mbox{ or } \zeta^2$ and
$c_2 = \zeta \mbox{ or } \zeta^2$, where $\zeta^3 = 1$.  Since $R_1
\neq 0$ when $(c_1,c_2) = (\zeta, \zeta^2)$, $(\zeta^2, \zeta)$, we
find that $(c_1,c_2) = (\zeta, \zeta) \mbox{ or } (\zeta^2,
\zeta^2)$.
\\
\indent Now assume $\textit{k}$ is not of characteristic 3,
$c_1^2-c_1+1=0$, and $R_1 = R_2 = 0$. Temporarily regard $c_1$
and $c_2$ as indeterminates, and temporarily
replace $R_1$ and $R_2$ with their natural preimages in $\mathbb{Z}[c_1,c_2]$. Using the Maple
command sprem, we recursively reduce the degree of
$c_1$ with initial polynomials $R_1$
and $c_1^2-c_1+1$, and $R_2$ and $c_1^2-c_1+1$. As a result: \\
$$(c_2^2-c_2+1)(c_2^2+c_2+1)(c_2^4-c_2^2+1)
(c_2^8-3c_2^7+9c_2^5+4c_2^4-18c_2^3+15c_2^2-6c_2+1)$$
$$(c_2^8-6c_2^7+15c_2^6-18c_2^5+4c_2^4+9c_2^3-3c_2+1) \in \langle R_1,
c_1^2-c_1+1 \rangle, \mbox{ and} \smallskip$$
$$c_2^2(3c_2^2-3c_2+1)(c_2^2-3c_2+3)(c_2^2-c_2+1)(c_2-1)^2 (4c_2^4+6c_2^3+c_2^2-3c_2+1)$$
$$(c_2^4-3c_2^3+c_2^2+6c_2+4) \in \langle R_2,
c_1^2-c_1+1 \rangle$$ \\
Note that if $R_1 = R_2 = 0$ and $c_1^2-c_1+1=0$ then
$c_2^2-c_2+1=0$. Now remove the temporary replacements of $c_1$,
$c_2$, $R_1$, and $R_2$. Thus if $R_1 = R_2 = 0$, then $c_1 = -\zeta$ or  $-\zeta^2$ and $c_2 = -\zeta$
or $-\zeta^2$, where $\zeta^3 = 1$.  However, if $(c_1,c_2) =
(-\zeta,-\zeta)$ or $(-\zeta^2,-\zeta^2)$, then $R_2 \neq 0$, which
is a contradiction.  Also, if $(c_1,c_2) = (-\zeta^2,-\zeta)$ or
$(-\zeta^2,-\zeta)$, then $R_1 \neq 0$, which is a
contradiction. \\
\indent We now consider all solutions to $R_1 = R_2= 0$.
\begin{itemize}
\begin{item}
$(c_1,c_2) = (1,1), (\zeta,\zeta), (\zeta^2,\zeta^2)$ if
$\textit{k}$ is of characteristic 2, \end{item}
\begin{item}
$(c_1,c_2) = (1,1),(1,-1), (-1,1), (-1,-1)$, if $\textit{k}$ is of
characteristic 3
\end{item}
\begin{item}
$(c_1,c_2) = (1,1), (1,-1), (-1,1), (-1,-1), (\zeta,\zeta),
(\zeta^2,\zeta^2)$ if $\textit{k}$ is not of characteristic 2 or 3.
\end{item}
\end{itemize}
Thus if $(X',Y')$ is not irreducible, $(c_1,c_2)$ is among the preceding pairs.
The pair of matrices $(X',Y')$ given when $(c_1,c_2) = (-1,-1),
(-1,1), \mbox{ and } (1,-1)$ are all equivalent. Also, the pair of
matrices $(X',Y')$ when $(c_1,c_2) = (\zeta,\zeta) \mbox{ and }
(\zeta^2,\zeta^2)$ are equivalent. Thus it is sufficient to check
irreducibility of $(X',Y')$ when $(c_1,c_2) = (1,1), (\zeta,\zeta),
\mbox{
and } (-1,-1).$  \\
\indent If $c_1=c_2$ and $c_1^3 = 1$, then $(X'Y')^2 = (Y'X')^2$.
Let $A = X'Y'$ and $B=Y'X'$. Then $(X',Y')$ is irreducible if and
only if $(A,B)$ is irreducible. Note $\Lambda \Gamma = \Gamma
\Lambda$ yields $A^6 = B^6 = I$ by Remark 3.3. Also, $ABA = BAB$.
From these relations, we can count possible monomials in $A$ and $B$
to find that $(A,B)$ span at most
an 18-dimensional space. Thus $(A,B)$ is not irreducible, and hence $(X',Y')$ is not irreducible. \\
If $c_1=c_2 = -1$ and $Q = \left( \begin{matrix} 1 & 1 & 0 & 0 & 0 &
0
\\ 0 & 0 & 0 & 0 & 1 & 1 \\ 0 & 0 & 1 & 1 & 0 & 0 \\ 1 & -1 & 0 & 0
& 0 & 0 \\ 0 & 0 & 0 & 0 & -1 & 1 \\ 0 & 0 & 1 & -1 & 0 & 0
\end{matrix} \right)$, then \bigskip \\ $$(QX'Q^{-1},QY'Q^{-1}) = \Biggl( \left(
\begin{matrix} 1 & & & & & \\ & 1 & & & & \\ & & 1 & & & \\ & & &
-1 & & \\ & & & & -1 & \\ & & & & & -1 \end{matrix} \right),
\left( \begin{matrix} 0 & 0 & 0  & 0 & 1 & 0 \\ 0 & 0 & 1 & 0 & 0
& 0 \\ 0 & 0 & 0 & -1 & 0 & 0 \\ 0 & -1 & 0 & 0 & 0 & 0 \\ 0 & 0 &
0 & 0 & 0 & -1 \\ -1 & 0 & 0 & 0 & 0 & 0 \end{matrix} \right)
\Biggr).$$ \\
Since $(QX'Q^{-1},QY'Q^{-1})$ is not irreducible, neither is
$(X',Y')$ when $c_1 = c_2 = -1$.  This gives the desired result.
\end{proof}
\begin{remark}
In Lemma 6.2, we proved further that if $$(X,Y) = \Biggl( \left(
\begin{matrix} 0 & 1 & 0 & 0 & 0 & 0 \\ 1 & 0 & 0 & 0 & 0 & 0 \\ 0
& 0 & 0 & 1 & 0 & 0 \\ 0 & 0 & 1 & 0 & 0 & 0 \\ 0 & 0 & 0 & 0 & 0
& 1 \\ 0 & 0 & 0 & 0 & 1 & 0 \end{matrix} \right), \left(
\begin{matrix} 0 & 0 & 0 & 0 & c_1 & 0 \\ 0 & 0 & 0 & 0 & 0 & 1 \\
c_2 & 0 & 0 & 0 & 0 & 0 \\ 0 & 1 & 0 & 0 & 0 & 0 \\ 0 & 0 &
\frac{1}{c_1c_2} & 0 & 0 & 0 \\ 0 & 0 & 0 & 1 & 0 & 0 \end{matrix}
\right) \Biggr)$$ where
\begin{enumerate}
\begin{item}
$(c_1,c_2) \neq (1,1), (-1,-1), (-1,1), (1,-1), (\zeta,\zeta),
(\zeta^2,\zeta^2)$, if $\textit{k}$ does not have characteristic 2
or 3, where $\zeta$ is a primitive cube root of unity;
\end{item}
\begin{item}
$(c_1,c_2) \neq (1,1), (\zeta,\zeta), (\zeta^2,\zeta^2)$, if
$\textit{k}$ has characteristic 2, where $\zeta$ is a primitive cube
root of unity;
\end{item}
\begin{item}
$(c_1,c_2) \neq (1,1), (-1,-1), (-1,1), (1,-1)$, if $\textit{k}$ has
characteristic 3,
\end{item}
\end{enumerate}
then $(X,Y)$ satisfies the hypotheses of Lemma 3.4.
\end{remark}
\noindent \textbf{6.4 } If $\rho: PSL_2(\mathbb{Z}) \rightarrow
GL_n(\textit{k})$ is an irreducible 6-dimensional representation of
$PSL_2(\mathbb{Z}) = \langle x,y$ $|$ $x^2 = y^3 = 1 \rangle$ with $X:=\rho(x)$ and
$Y:=\rho(y)$, it follows from Lemma 6.1 that $X$ and $Y$ have one nonzero entry per
row and column.  In particular, $\rho: PSL_2(\mathbb{Z}) \rightarrow
GL_n(\textit{k})$ is an irreducible 6-dimensional representation of
$PSL_2(\mathbb{Z})$ if and only if $(X,Y)$ satisfies the hypothesis of Lemma 3.4.  Then
by Lemma 6.2 and Remark 6.3, we have classified the irreducible 6-dimensional representations of $PSL_2(\mathbb{Z})$
up to equivalence.

\section{Appendix 1}
In this Appendix, we list the matrices used in the case-by-case
analysis proofs in Lemma 6.1. Again, let $\gamma_i$ denote
$\Gamma_{i,i},$ and $\lambda_i$ denote $\Lambda_{i,i}$.
Let $(a_1, a_2, a_3, a_4, a_5, a_6)$ denote the
unique $6 \times 6$ permutation matrix which maps $e_i$ to
$e_{a_i}$, where $e_i$ is the
standard basis vector of $\textit{k}^6$. \\
\textbf{Case 5:} In our proof, we found that $\lambda_1 =
\frac{1}{\lambda_1}$ and $\lambda_5 = \frac{1}{\lambda_6}$.  From $\Gamma X \Gamma = X$,
we conclude that if $\gamma_i \gamma_j \neq 1$, then $X_{i,j} = 0$.  From $\Lambda X \Lambda = X$,
we conclude that if $\lambda_i \lambda_j \neq 1$, then $X_{i,j} = 0$. For each of the
30 matrices for $\Gamma$ in terms of $\lambda_i$, we determine which entries of $X$
must equal 0.  Of the 30 possible matrices for $\Gamma$, we found that for all but 14 matrices, $X$
is forced to be singular.  Two of these $\Gamma$ lead to $(X,Y)$ irreducible.  For
the 12 remaining matrices for $\Gamma$, we replace $(X,Y)$ with $(PXP^{-1},PYP^{-1})$
for an appropriate permutation matrix $P$, so that $\Lambda$ is preserved, $\gamma_1 =$
$\gamma_2 =$ $\gamma_5 =$ $\gamma_6 = \lambda_1$, $\gamma_3 =
\lambda_5$, and $\gamma_4 = \lambda_6$.  We list each of the 14 possible matrices for $\Gamma$, and
for the 12 possible matrices for $\Gamma$ which do not lead to $(X,Y)$ irreducible, we list the
corresponding permutation matrix $P$.
\begin{enumerate}
\begin{item}
$\gamma_1 = \gamma_2 = \gamma_3 = \gamma_4 = \lambda_1$, $\gamma_5 =
\lambda_5$, $\gamma_6 = \lambda_6$; $\Lambda = \Gamma$ implies
$(X,Y)$ reducible.
\end{item}
\begin{item}
$\gamma_1 = \gamma_2 = \gamma_3 = \gamma_4 = \lambda_1$, $\gamma_5 =
\lambda_6$, $\gamma_6 = \lambda_5$; $\Lambda \Gamma = I$ implies
$(X,Y)$ reducible.
\end{item}
\begin{item}
$\gamma_1 = \gamma_2 = \gamma_5 = \gamma_6 = \lambda_1$, $\gamma_3 =
\lambda_5$, $\gamma_4 = \lambda_6$; $P$ is the identity matrix.
\end{item}
\begin{item}
$\gamma_1 = \gamma_2 = \gamma_5 = \gamma_6 = \lambda_1$, $\gamma_3 =
\lambda_6$, $\gamma_4 = \lambda_5$; $P=(1, 2, 4, 3, 5, 6)$.
\end{item}
\begin{item}
$\gamma_1 = \gamma_3 = \gamma_5 = \gamma_6 = \lambda_1$, $\gamma_2 =
\lambda_5$, $\gamma_4 = \lambda_6$; $P=(1, 3, 2, 4, 5, 6)$.
\end{item}
\begin{item}
$\gamma_1 = \gamma_4 = \gamma_5 = \gamma_6 = \lambda_1$, $\gamma_2 =
\lambda_5$, $\gamma_3 = \lambda_6$; $P=(1, 4, 2, 3, 5, 6)$.
\end{item}
\begin{item}
$\gamma_1 = \gamma_4 = \gamma_5 = \gamma_6 = \lambda_1$, $\gamma_2 =
\lambda_6$, $\gamma_3 = \lambda_5$; $P=(1, 4, 3, 2, 5, 6)$.
\end{item}
\begin{item}
$\gamma_1 = \gamma_3 = \gamma_5 = \gamma_6 = \lambda_1$, $\gamma_2 =
\lambda_6$, $\gamma_4 = \lambda_5$; $P=(1, 3, 4, 2, 5, 6)$.
\end{item}
\begin{item}
$\gamma_2 = \gamma_3 = \gamma_5 = \gamma_6 = \lambda_1$, $\gamma_1 =
\lambda_5$, $\gamma_4 = \lambda_6$; $P=(2, 3, 1, 4, 5, 6)$.
\end{item}
\begin{item}
$\gamma_2 = \gamma_4 = \gamma_5 = \gamma_6 = \lambda_1$, $\gamma_1 =
\lambda_5$, $\gamma_3 = \lambda_6$; $P=(2, 4, 1, 3, 5, 6)$.
\end{item}
\begin{item}
$\gamma_3 = \gamma_4 = \gamma_5 = \gamma_6 = \lambda_1$, $\gamma_1 =
\lambda_5$, $\gamma_2 = \lambda_6$; $P=(3, 4, 1, 2, 5, 6)$.
\end{item}
\begin{item}
$\gamma_2 = \gamma_3 = \gamma_5 = \gamma_6 = \lambda_1$, $\gamma_1 =
\lambda_6$, $\gamma_4 = \lambda_5$; $P=(2, 3, 4, 1, 5, 6)$.
\end{item}
\begin{item}
$\gamma_2 = \gamma_4 = \gamma_5 = \gamma_6 = \lambda_1$, $\gamma_1 =
\lambda_6$, $\gamma_3 = \lambda_5$; $P=(2, 4, 3, 1, 5, 6)$.
\end{item}
\begin{item}
$\gamma_3 = \gamma_4 = \gamma_5 = \gamma_6 = \lambda_1$, $\gamma_1 =
\lambda_6$, $\gamma_2 = \lambda_5$; $P=(3, 4, 2, 1, 5, 6)$.
\end{item}
\end{enumerate}
\textbf{Case 7:} Here we list the 22 possible values for $\Gamma$
and the corresponding forms of $X$ and $Y$.  The forms of $X$ and
$Y$ are described by which entries (of $X$ or $Y$ respectively) must
equal 0.  There are 5 forms for $X$ and 22 forms for $Y$. \\
$(\gamma_1,\gamma_2,\gamma_3,\gamma_4,\gamma_5,\gamma_6) = $ \\
$\begin{matrix} (1) \mbox{ }
(\lambda_1,\lambda_1,\lambda_3,\lambda_3,\lambda_5,\lambda_5) & (2)
\mbox{ }
(\lambda_1,\lambda_1,\lambda_3,\lambda_5,\lambda_5,\lambda_3) & (3)
\mbox{ }
(\lambda_1,\lambda_1,\lambda_3,\lambda_5,\lambda_3,\lambda_5) \\ (4)
\mbox{ }
(\lambda_1,\lambda_1,\lambda_5,\lambda_3,\lambda_3,\lambda_5) & (5)
\mbox{ }
(\lambda_1,\lambda_1,\lambda_5,\lambda_3,\lambda_5,\lambda_3) & (6)
\mbox{ }
(\lambda_1,\lambda_1,\lambda_5,\lambda_5,\lambda_3,\lambda_3) \\ (7)
\mbox{ }
(\lambda_3,\lambda_5,\lambda_3,\lambda_1,\lambda_5,\lambda_1) & (8)
\mbox{ }
(\lambda_3,\lambda_5,\lambda_3,\lambda_1,\lambda_1,\lambda_5) & (9)
\mbox{ }
(\lambda_3,\lambda_5,\lambda_5,\lambda_1,\lambda_3,\lambda_1) \\
(10) \mbox{ }
(\lambda_3,\lambda_5,\lambda_5,\lambda_1,\lambda_1,\lambda_3) & (11)
\mbox{ }
(\lambda_3,\lambda_5,\lambda_1,\lambda_3,\lambda_5,\lambda_1) & (12)
\mbox{ }
(\lambda_3,\lambda_5,\lambda_1,\lambda_3,\lambda_1,\lambda_5) \\
(13) \mbox{ }
(\lambda_3,\lambda_5,\lambda_1,\lambda_5,\lambda_3,\lambda_1) & (14)
\mbox{ }
(\lambda_3,\lambda_5,\lambda_1,\lambda_5,\lambda_1,\lambda_3) & (15)
\mbox{ }
(\lambda_5,\lambda_3,\lambda_1,\lambda_3,\lambda_1,\lambda_5) \\
(16) \mbox{ }
(\lambda_5,\lambda_3,\lambda_1,\lambda_3,\lambda_5,\lambda_1) & (17)
\mbox{ }
(\lambda_5,\lambda_3,\lambda_1,\lambda_5,\lambda_1,\lambda_3) & (18)
\mbox{ }
(\lambda_5,\lambda_3,\lambda_1,\lambda_5,\lambda_3,\lambda_1) \\
(19) \mbox{ }
(\lambda_5,\lambda_3,\lambda_3,\lambda_1,\lambda_1,\lambda_5) & (20)
\mbox{ }
(\lambda_5,\lambda_3,\lambda_3,\lambda_1,\lambda_5,\lambda_1) & (21)
\mbox{ }
(\lambda_5,\lambda_3,\lambda_5,\lambda_1,\lambda_1,\lambda_3) \\
(22) \mbox{ }
(\lambda_5,\lambda_3,\lambda_5,\lambda_1,\lambda_3,\lambda_1) & &
\end{matrix}$ \\ \\
For $\Gamma$ in cases 1 and 6, $X$ has the following form: if
$X_{i,j} \neq 0$, then $(i,j) = (1,1),$ $(1,2),$ $(2,1),$ $(2,2),$
$(3,5),$ $(3,6),$ $(4,5),$ $(4,6),$ $(5,3),$ $(5,4),$ $(6,3),$ or
$(6,4)$. \\ For $\Gamma$ in cases 2 and 4, $X$ has the following
form: if $X_{i,j} \neq 0$, then $(i,j) = (1,1),$ $(1,2),$ $(2,1),$
$(2,2),$ $(3,5),$ $(4,6),$ $(5,3),$ or $(6,4)$. \\ For $\Gamma$ in
cases 3 and 5, $X$ has the following form: if $X_{i,j} \neq 0$, then
$(i,j) = (1,1),$ $(1,2),$ $(2,1),$ $(2,2),$ $(3,6),$ $(4,5),$
$(5,4),$ or $(6,3)$.  \\ For $\Gamma$ in cases 7, 9, 12, 14, 15, 17,
20, and 22, $X$ has the following form: if $X_{i,j} \neq 0$, then
$(i,j) = (1,2),$ $(2,1),$ $(3,5),$ $(4,6),$ $(5,3),$ or $(6,4)$.
\\ For $\Gamma$ in cases 8, 10, 11, 13, 16, 18, 19, and 21, $X$ has
the following form: if $X_{i,j} \neq 0$, then $(i,j) = (1,2),$
$(2,1),$ $(3,6),$ $(4,5),$ $(5,4),$ or $(6,3)$. \\ Fix $\Gamma$ in
one of the above cases.  Then $Y$ has the following form: if
$Y_{i,j} \neq 0$, then $\lambda_i \gamma_j = 1$.  Note that
$\lambda_1 = \frac{1}{\lambda_1}$ and $\lambda_3 =
\frac{1}{\lambda_5}$. \\  The selection of $\Gamma$ in the cases 7, 8, 9, 10, 11,
12, 13, 14, 15, 16, 17, 18, 19, 20, 21, and 22 does not force $(X,Y)$ to be reducible.

\textbf{Case 8:} Below are the 20 possible values for $\Gamma$ and
the corresponding forms of $X$ and $Y$.  The forms of $X$ and $Y$
are described by which entries (of $X$ or $Y$ respectively) must
equal 0.  There are 7 forms for $X$ and 20 forms for $Y$. \\
$(\gamma_1,\gamma_2,\gamma_3,\gamma_4,\gamma_5,\gamma_6) = $ \\
$\begin{matrix} (1) \mbox{ }
(\lambda_1,\lambda_5,\lambda_6,\lambda_4,\lambda_1,\lambda_1) & (2)
\mbox{ }
(\lambda_1,\lambda_6,\lambda_5,\lambda_4,\lambda_1,\lambda_1) & (3)
\mbox{ }
(\lambda_5,\lambda_1,\lambda_6,\lambda_4,\lambda_1,\lambda_1) \\ (4)
\mbox{ }
(\lambda_5,\lambda_6,\lambda_1,\lambda_4,\lambda_1,\lambda_1) & (5)
\mbox{ }
(\lambda_6,\lambda_1,\lambda_5,\lambda_4,\lambda_1,\lambda_1) & (6)
\mbox{ }
(\lambda_6,\lambda_5,\lambda_1,\lambda_4,\lambda_1,\lambda_1) \\ (7)
\mbox{ }
(\lambda_4,\lambda_5,\lambda_6,\lambda_1,\lambda_1,\lambda_1) & (8)
\mbox{ }
(\lambda_4,\lambda_6,\lambda_5,\lambda_1,\lambda_1,\lambda_1) & (9)
\mbox{ }
(\lambda_5,\lambda_4,\lambda_6,\lambda_1,\lambda_1,\lambda_1) \\
(10) \mbox{ }
(\lambda_5,\lambda_6,\lambda_4,\lambda_1,\lambda_1,\lambda_1) & (11)
\mbox{ }
(\lambda_6,\lambda_4,\lambda_5,\lambda_1,\lambda_1,\lambda_1) & (12)
\mbox{ }
(\lambda_6,\lambda_5,\lambda_4,\lambda_1,\lambda_1,\lambda_1) \\
(13) \mbox{ }
(\lambda_1,\lambda_1,\lambda_1,\lambda_4,\lambda_5,\lambda_6) & (14)
\mbox{ }
(\lambda_1,\lambda_1,\lambda_4,\lambda_1,\lambda_5,\lambda_6) & (15)
\mbox{ }
(\lambda_1,\lambda_4,\lambda_1,\lambda_1,\lambda_5,\lambda_6) \\
(16) \mbox{ }
(\lambda_4,\lambda_1,\lambda_1,\lambda_1,\lambda_5,\lambda_6) & (17)
\mbox{ }
(\lambda_1,\lambda_1,\lambda_1,\lambda_4,\lambda_6,\lambda_5) & (18)
\mbox{ }
(\lambda_1,\lambda_1,\lambda_4,\lambda_1,\lambda_6,\lambda_5) \\
(19) \mbox{ }
(\lambda_1,\lambda_4,\lambda_1,\lambda_1,\lambda_6,\lambda_5) & (20)
\mbox{ }
(\lambda_4,\lambda_1,\lambda_1,\lambda_1,\lambda_6,\lambda_5) &
\end{matrix}$ \\ \\
For $\Gamma$ in cases 1, 2, 7, and 8, $X$ has the following form: if
$X_{i,j} \neq 0$, then $(i,j) = (1,1),$ $(2,3),$ $(3,2),$ $(4,4),$
$(5,6),$ or $(6,5)$. \\
For $\Gamma$ in cases 3, 5, 9, and 11, $X$ has the following form:
if $X_{i,j} \neq 0$, then $(i,j) = (1,3),$ $(2,2),$ $(3,1),$
$(4,4),$ $(5,6),$ or $(6,5)$. \\ For $\Gamma$ in cases 4, 6, 10, and
12, $X$ has the following form: if $X_{i,j} \neq 0$, then $(i,j) =
(1,2),$ $(2,1),$ $(3,3),$ $(4,4),$ $(5,6),$ or $(6,5)$. \\
For $\Gamma$ in cases 13 and 17, $X$ has the following form: if
$X_{i,j} \neq 0$, then $(i,j) = (1,1),$ $(1,2),$ $(1,3),$ $(2,1),$
$(2,2),$ $(2,3)$, $(3,1),$ $(3,2),$ $(3,3),$ $(4,4),$ $(5,6),$ or
$(6,5)$.
\\ For $\Gamma$ in cases 14 and 18, $X$ has
the following form: if $X_{i,j} \neq 0$, then $(i,j) = (1,1),$
$(1,2),$ $(2,1),$ $(2,2),$ $(3,3),$ $(4,4),$ $(5,6)$, or $(6,5)$.
\\ For $\Gamma$ in cases 15 and 19, $X$ has the following form: if
$X_{i,j} \neq 0$, then $(i,j) = (1,1),$ $(1,3),$ $(2,2),$ $(3,1),$
$(3,3),$ $(4,4),$ $(5,6)$, or $(6,5)$.
\\ For $\Gamma$ in cases 16 and 20, $X$ has the following form: if
$X_{i,j} \neq 0$, then $(i,j) = (1,1),$
$(2,2),$ $(2,3),$ $(3,2),$ $(3,3),$ $(4,4),$ $(5,6)$, or $(6,5)$. \\
Fix $\Gamma$ in one of the above cases.  Then $Y$ has the following
form: if $Y_{i,j} \neq 0$, then $\lambda_i \gamma_j = 1$.  Note that
$\lambda_1 = \frac{1}{\lambda_1}$, $\lambda_4 =
\frac{1}{\lambda_4}$, and $\lambda_5 =
\frac{1}{\lambda_6}$. \\
The selection of $\Gamma$ in the cases 7, 8, 9, 10, 11, and 12, does not force $(X,Y)$ to be reducible. \\

\textbf{Case 9:} Here, we list the possible values for $\Gamma$ and
the corresponding forms of $X$ and $Y$ for two values of $\Lambda$.
\\ First we assume $\lambda_1 = \frac{1}{\lambda_1}$, $\lambda_3 =
\frac{1}{\lambda_3}$, and $\lambda_5 = \frac{1}{\lambda_6}$. There
are 20 possible values of $\Gamma$.  The forms of $X$ and $Y$ are
described by which entries (of $X$ or $Y$ respectively) must equal
0.  There are 4 forms for $X$ and 20 forms for $Y$. \\
$(\gamma_1,\gamma_2,\gamma_3,\gamma_4,\gamma_5,\gamma_6) = $ \\
$\begin{matrix} (1) \mbox{ }
(\lambda_3,\lambda_3,\lambda_5,\lambda_6,\lambda_1,\lambda_1) & (2)
\mbox{ }
(\lambda_3,\lambda_3,\lambda_6,\lambda_5,\lambda_1,\lambda_1) & (3)
\mbox{ }
(\lambda_5,\lambda_6,\lambda_3,\lambda_3,\lambda_1,\lambda_1) \\ (4)
\mbox{ }
(\lambda_6,\lambda_5,\lambda_3,\lambda_3,\lambda_1,\lambda_1) & (5)
\mbox{ }
(\lambda_1,\lambda_1,\lambda_5,\lambda_6,\lambda_3,\lambda_3) & (6)
\mbox{ }
(\lambda_1,\lambda_1,\lambda_6,\lambda_5,\lambda_3,\lambda_3) \\ (7)
\mbox{ }
(\lambda_5,\lambda_6,\lambda_1,\lambda_1,\lambda_3,\lambda_3) & (8)
\mbox{ }
(\lambda_6,\lambda_5,\lambda_1,\lambda_1,\lambda_3,\lambda_3) & (9)
\mbox{ }
(\lambda_1,\lambda_1,\lambda_3,\lambda_3,\lambda_5,\lambda_6) \\
(10) \mbox{ }
(\lambda_1,\lambda_1,\lambda_3,\lambda_3,\lambda_6,\lambda_5) & (11)
\mbox{ }
(\lambda_1,\lambda_3,\lambda_1,\lambda_3,\lambda_5,\lambda_6) & (12)
\mbox{ }
(\lambda_1,\lambda_3,\lambda_1,\lambda_3,\lambda_6,\lambda_5) \\
(13) \mbox{ }
(\lambda_1,\lambda_3,\lambda_3,\lambda_1,\lambda_5,\lambda_6) & (14)
\mbox{ }
(\lambda_1,\lambda_3,\lambda_3,\lambda_1,\lambda_6,\lambda_5) & (15)
\mbox{ }
(\lambda_3,\lambda_1,\lambda_1,\lambda_3,\lambda_5,\lambda_6) \\
(16) \mbox{ }
(\lambda_3,\lambda_1,\lambda_1,\lambda_3,\lambda_6,\lambda_5) & (17)
\mbox{ }
(\lambda_3,\lambda_1,\lambda_3,\lambda_1,\lambda_5,\lambda_6) & (18)
\mbox{ }
(\lambda_3,\lambda_1,\lambda_3,\lambda_1,\lambda_6,\lambda_5) \\
(19) \mbox{ }
(\lambda_3,\lambda_3,\lambda_1,\lambda_1,\lambda_5,\lambda_6) & (20)
\mbox{ }
(\lambda_3,\lambda_3,\lambda_1,\lambda_1,\lambda_6,\lambda_5) &
\end{matrix}$ \\ \\
For $\Gamma$ in cases 3, 4, 7, and 8, $X$ has the following form: if
$X_{i,j} \neq 0$, then $(i,j) = (1,2),$ $(2,1),$ $(3,3),$ $(3,4),$
$(4,3),$ $(4,4),$ $(5,6),$ or $(6,5)$. \\
For $\Gamma$ in cases 1, 2, 5, and 6, $X$ has the following form: if
$X_{i,j} \neq 0$, then $(i,j) = (1,1),$ $(1,2),$ $(2,1),$ $(2,2),$
$(3,4),$ $(4,3),$ $(5,6),$ or $(6,5)$. \\
For $\Gamma$ in cases 9, 10, 19, and 20, $X$ has the following form:
if $X_{i,j} \neq 0$, then $(i,j) = (1,1),$ $(1,2),$ $(2,1),$
$(2,2),$ $(3,3),$ $(3,4),$ $(4,3),$ $(4,4),$ $(5,6),$ or $(6,5)$. \\
For $\Gamma$ in cases 11, 12, 13, 14, 15, 16, 17, and 18, $X$ has
the following form: if $X_{i,j} \neq 0$, then $(i,j) = (1,1),$
$(2,2),$ $(3,3),$ $(4,4),$ $(5,6),$ or $(6,5)$. \\ Fix $\Gamma$ in
one of the above cases.  Then $Y$ has the following form: if
$Y_{i,j} \neq 0$, then $\lambda_i \gamma_j = 1$.
Note that the selection of $\Gamma$ in the cases 1, 2, 7, and 8, does
not force $(X,Y)$ to be reducible. \\

Now assume $\lambda_1 = \frac{1}{\lambda_3}$ and $\lambda_5 =
\frac{1}{\lambda_6}$. There are 44 possible values of $\Gamma$. The
forms of $X$ and $Y$ are described by which entries (of $X$ or $Y$
respectively) must equal
0.  There are 3 forms for $X$ and 44 forms for $Y$. \\
$(\gamma_1,\gamma_2,\gamma_3,\gamma_4,\gamma_5,\gamma_6) = $ \\
$\begin{matrix} (1) \mbox{ }
(\lambda_1,\lambda_5,\lambda_3,\lambda_6,\lambda_1,\lambda_3) & (2)
\mbox{ }
(\lambda_1,\lambda_5,\lambda_6,\lambda_3,\lambda_1,\lambda_3) & (3)
\mbox{ }
(\lambda_1,\lambda_6,\lambda_3,\lambda_5,\lambda_1,\lambda_3) \\ (4)
\mbox{ }
(\lambda_1,\lambda_6,\lambda_5,\lambda_3,\lambda_1,\lambda_3) & (5)
\mbox{ }
(\lambda_3,\lambda_5,\lambda_1,\lambda_6,\lambda_1,\lambda_3) & (6)
\mbox{ }
(\lambda_3,\lambda_5,\lambda_6,\lambda_1,\lambda_1,\lambda_3) \\ (7)
\mbox{ }
(\lambda_3,\lambda_6,\lambda_1,\lambda_5,\lambda_1,\lambda_3) & (8)
\mbox{ }
(\lambda_3,\lambda_6,\lambda_5,\lambda_1,\lambda_1,\lambda_3) & (9)
\mbox{ }
(\lambda_5,\lambda_1,\lambda_3,\lambda_6,\lambda_1,\lambda_3) \\
(10) \mbox{ }
(\lambda_5,\lambda_1,\lambda_6,\lambda_3,\lambda_1,\lambda_3) & (11)
\mbox{ }
(\lambda_5,\lambda_3,\lambda_1,\lambda_6,\lambda_1,\lambda_3) & (12)
\mbox{ }
(\lambda_5,\lambda_3,\lambda_6,\lambda_1,\lambda_1,\lambda_3) \\
(13) \mbox{ }
(\lambda_6,\lambda_1,\lambda_3,\lambda_5,\lambda_1,\lambda_3) & (14)
\mbox{ }
(\lambda_6,\lambda_1,\lambda_5,\lambda_3,\lambda_1,\lambda_3) & (15)
\mbox{ }
(\lambda_6,\lambda_3,\lambda_1,\lambda_5,\lambda_1,\lambda_3) \\
(16) \mbox{ }
(\lambda_6,\lambda_3,\lambda_5,\lambda_1,\lambda_1,\lambda_3) & (17)
\mbox{ }
(\lambda_1,\lambda_5,\lambda_3,\lambda_6,\lambda_3,\lambda_1) & (18)
\mbox{ }
(\lambda_1,\lambda_5,\lambda_6,\lambda_3,\lambda_3,\lambda_1) \\
(19) \mbox{ }
(\lambda_1,\lambda_6,\lambda_3,\lambda_5,\lambda_3,\lambda_1) & (20)
\mbox{ }
(\lambda_1,\lambda_6,\lambda_5,\lambda_3,\lambda_3,\lambda_1) & (21)
\mbox{ }
(\lambda_3,\lambda_5,\lambda_1,\lambda_6,\lambda_3,\lambda_1) \\
(22) \mbox{ }
(\lambda_3,\lambda_5,\lambda_6,\lambda_1,\lambda_3,\lambda_1) & (23)
\mbox{ }
(\lambda_3,\lambda_6,\lambda_1,\lambda_5,\lambda_3,\lambda_1) & (24)
\mbox{ }
(\lambda_3,\lambda_6,\lambda_5,\lambda_1,\lambda_3,\lambda_1) \\
(25) \mbox{ }
(\lambda_5,\lambda_1,\lambda_3,\lambda_6,\lambda_3,\lambda_1) & (26)
\mbox{ }
(\lambda_5,\lambda_1,\lambda_6,\lambda_3,\lambda_3,\lambda_1) & (27)
\mbox{ }
(\lambda_5,\lambda_3,\lambda_1,\lambda_6,\lambda_3,\lambda_1) \\
(28) \mbox{ }
(\lambda_5,\lambda_3,\lambda_6,\lambda_1,\lambda_3,\lambda_1) & (29)
\mbox{ }
(\lambda_6,\lambda_1,\lambda_3,\lambda_5,\lambda_3,\lambda_1) & (30)
\mbox{ }
(\lambda_6,\lambda_1,\lambda_5,\lambda_3,\lambda_3,\lambda_1) \\
(31) \mbox{ }
(\lambda_6,\lambda_3,\lambda_1,\lambda_5,\lambda_3,\lambda_1) & (32)
\mbox{ }
(\lambda_6,\lambda_3,\lambda_5,\lambda_1,\lambda_3,\lambda_1) & (33)
\mbox{ }
(\lambda_1,\lambda_1,\lambda_3,\lambda_3,\lambda_5,\lambda_6) \\
(34) \mbox{ }
(\lambda_1,\lambda_3,\lambda_1,\lambda_3,\lambda_5,\lambda_6) & (35)
\mbox{ }
(\lambda_1,\lambda_3,\lambda_3,\lambda_1,\lambda_5,\lambda_6) & (36)
\mbox{ }
(\lambda_3,\lambda_1,\lambda_1,\lambda_3,\lambda_5,\lambda_6) \\
(37) \mbox{ }
(\lambda_3,\lambda_1,\lambda_3,\lambda_1,\lambda_5,\lambda_6) & (38)
\mbox{ }
(\lambda_3,\lambda_3,\lambda_1,\lambda_1,\lambda_5,\lambda_6) & (39)
\mbox{ }
(\lambda_1,\lambda_1,\lambda_3,\lambda_3,\lambda_6,\lambda_5) \\
(40) \mbox{ }
(\lambda_1,\lambda_3,\lambda_1,\lambda_3,\lambda_6,\lambda_5) & (41)
\mbox{ }
(\lambda_1,\lambda_3,\lambda_3,\lambda_1,\lambda_6,\lambda_5) & (42)
\mbox{ }
(\lambda_3,\lambda_1,\lambda_1,\lambda_3,\lambda_6,\lambda_5) \\
(43) \mbox{ }
(\lambda_3,\lambda_1,\lambda_3,\lambda_1,\lambda_6,\lambda_5) & (44)
\mbox{ }
(\lambda_3,\lambda_3,\lambda_1,\lambda_1,\lambda_6,\lambda_5) &
\end{matrix}$ \\ \\
For $\Gamma$ in cases 1, 3, 5, 7, 10, 12, 14, 16, 17, 19, 21, 23, 26, 28, 30, 32, 35, 36, 41, and 42,
$X$ has the following form: if $X_{i,j} \neq 0$, then $(i,j) = (1,3),$ $(2,4),$ $(3,1),$ $(4,2),$ $(5,6),$ or $(6,5)$. \\
For $\Gamma$ in cases 2, 4, 6, 8, 9, 11, 13, 15, 18, 20, 22, 24, 25, 27, 29, 31, 34, 37, 40, and 43
$X$ has the following form: if $X_{i,j} \neq 0$, then $(i,j) = (1,4),$ $(2,3),$ $(3,2),$ $(4,1),$ $(5,6),$ or $(6,5)$. \\
For $\Gamma$ in cases 33, 38, 39, and 44, $X$ has the following form:
if $X_{i,j} \neq 0$, then $(i,j) = (1,3),$ $(1,4),$ $(2,3),$
$(2,4),$ $(3,1),$ $(3,2),$ $(4,1),$ $(4,2),$ $(5,6),$ or $(6,5)$. \\
Fix $\Gamma$ in one of the above cases.  Then $Y$ has the following form: if
$Y_{i,j} \neq 0$,
then $\lambda_i \gamma_j = 1$. We note that the selection of $\Gamma$ in the
cases 1 through 32 does not force $(X,Y)$ to be reducible. \\

\textbf{Case 10:} Listed below are the 24 possible values for $\Gamma$ and
the corresponding forms of $X$ and $Y$.  The forms of $X$ and $Y$
are described by which entries (of $X$ or $Y$ respectively) must
equal 0.  There are 2 forms for $X$ and 24 forms for $Y$. \\
$(\gamma_1,\gamma_2,\gamma_3,\gamma_4,\gamma_5,\gamma_6) = $ \\
$\begin{matrix} (1) \mbox{ }
(\lambda_1,\lambda_1,\lambda_5,\lambda_6,\lambda_3,\lambda_4) & (2)
\mbox{ }
(\lambda_1,\lambda_1,\lambda_6,\lambda_5,\lambda_3,\lambda_4) & (3)
\mbox{ }
(\lambda_5,\lambda_6,\lambda_1,\lambda_1,\lambda_3,\lambda_4) \\ (4)
\mbox{ }
(\lambda_6,\lambda_5,\lambda_1,\lambda_1,\lambda_3,\lambda_4) & (5)
\mbox{ }
(\lambda_1,\lambda_1,\lambda_5,\lambda_6,\lambda_4,\lambda_3) & (6)
\mbox{ }
(\lambda_1,\lambda_1,\lambda_6,\lambda_5,\lambda_4,\lambda_3) \\ (7)
\mbox{ }
(\lambda_5,\lambda_6,\lambda_1,\lambda_1,\lambda_4,\lambda_3) & (8)
\mbox{ }
(\lambda_6,\lambda_5,\lambda_1,\lambda_1,\lambda_4,\lambda_3) & (9)
\mbox{ }
(\lambda_1,\lambda_1,\lambda_3,\lambda_4,\lambda_5,\lambda_6) \\
(10) \mbox{ }
(\lambda_1,\lambda_1,\lambda_4,\lambda_3,\lambda_5,\lambda_6) & (11)
\mbox{ }
(\lambda_3,\lambda_4,\lambda_1,\lambda_1,\lambda_5,\lambda_6) & (12)
\mbox{ }
(\lambda_4,\lambda_3,\lambda_1,\lambda_1,\lambda_5,\lambda_6) \\
(13) \mbox{ }
(\lambda_1,\lambda_1,\lambda_3,\lambda_4,\lambda_6,\lambda_5) & (14)
\mbox{ }
(\lambda_1,\lambda_1,\lambda_4,\lambda_3,\lambda_6,\lambda_5) & (15)
\mbox{ }
(\lambda_3,\lambda_4,\lambda_1,\lambda_1,\lambda_6,\lambda_5) \\
(16) \mbox{ }
(\lambda_4,\lambda_3,\lambda_1,\lambda_1,\lambda_6,\lambda_5) & (17)
\mbox{ }
(\lambda_3,\lambda_4,\lambda_5,\lambda_6,\lambda_1,\lambda_1) & (18)
\mbox{ }
(\lambda_3,\lambda_4,\lambda_6,\lambda_5,\lambda_1,\lambda_1) \\
(19) \mbox{ }
(\lambda_4,\lambda_3,\lambda_5,\lambda_6,\lambda_1,\lambda_1) & (20)
\mbox{ }
(\lambda_4,\lambda_3,\lambda_6,\lambda_5,\lambda_1,\lambda_1) & (21)
\mbox{ }
(\lambda_5,\lambda_6,\lambda_3,\lambda_4,\lambda_1,\lambda_1) \\ (22)
\mbox{ }
(\lambda_5,\lambda_6,\lambda_4,\lambda_3,\lambda_1,\lambda_1) & (23)
\mbox{ }
(\lambda_6,\lambda_5,\lambda_3,\lambda_4,\lambda_1,\lambda_1) & (24)
\mbox{ }
(\lambda_6,\lambda_5,\lambda_4,\lambda_3,\lambda_1,\lambda_1)
\end{matrix}$ \\ \\
For $\Gamma$ in cases 3, 4, 7, 8, 11, 12, 15, 16, 17, 18, 19, 20, 21, 22, 23, and 24,
$X$ has the following form: if $X_{i,j} \neq 0$, then $(i,j) = (1,2),$ $(2,1),$ $(3,4),$ $(4,3),$
$(5,6),$ or $(6,5)$. \\
For $\Gamma$ in cases 1, 2, 5, 6, 9, 10, 13, and 14, $X$ has the following form:
if $X_{i,j} \neq 0$, then $(i,j) = (1,1),$ $(1,2),$ $(2,1),$
$(2,2),$ $(3,4),$ $(4,3),$ $(5,6),$ or $(6,5)$. \\
Fix $\Gamma$ in one of the above cases.  Then $Y$ has the following
form: if $Y_{i,j} \neq 0$, then $\lambda_i \gamma_j = 1$.  Note that
$\lambda_1 = \frac{1}{\lambda_1}$, $\lambda_3 =
\frac{1}{\lambda_4}$, and $\lambda_5 =
\frac{1}{\lambda_6}$.  Also, the selection of $\Gamma$ in the cases 3, 4, 7, 8, 17, 18,
19, and 20, does not force $(X,Y)$ to be reducible. \\

\section{Appendix 2}
\noindent In this Appendix, we list the specific choices of polynomials
and the Maple commands used in the latter half of the proof of Lemma 6.2.
\\ Our choices of $F_1$, $F_2$, and $F_3$ are:
$$F_1=c_1^2c_2^2 - c_1^2c_2 +c_1^2 -c_1c_2 -c_1 +1$$
$$F_2= c_1^2c_2^2 - c_1c_2^2 + c_2^2 - c_1c_2 - c_2 +1 \smallskip$$
$$F_3= c_1^4c_2^4 +c_1^4c_2^3 +c_1^4c_2^2 +c_1^3c_2^4 -c_1^3c_2^3
-c_1^3c_2^2 +c_1^3c_2 +c_1^2c_2^4 -c_1^2c_2^2-6c_1^2c_2^2$$
$$-c_1^2c_2 +c_1^2 +c_1c_2^3 -c_1c_2^2 -c_1c_2 +c_1
+c_2^2 +c_2 +1.$$

For multivariate polynomials $A$ and $B$ in the indeterminate $x$, sprem($A$, $B$,
$x$, `m', `q') outputs a multivariate polynomial $r$, where $mA = qB
+ r$ and $m, q, r$ are polynomials over $\mathbb{Z}$. The variables
$m$ and $q$ are assigned their corresponding polynomial values.
Furthermore, the degree of $x$ in $r$ is strictly less than the
degree of $x$ in $B$. We use the command `factor' to determine the
roots of the polynomials in question. \\
\\
To solve $R_1=F_1=0$: \\
{\tt $>$ A1:= sprem(R1, F1, c2, `m', `q'); \\
$>$ sprem(F1, A1, c2, `m', `q');\\
$>$ factor(\%);}
\\
\\
To solve $R_1=F_2=0$:\\
{\tt $>$ A2:= sprem(R1, F2, c2, `m', `q');\\
$>$ sprem(F2, A2, c2, `m', `q');\\
$>$ factor(\%);\\
\\}
To solve $R_1=F_3=0$: \\
\indent Without factoring $R_1$:\\
{\tt $>$ A3:= sprem(R1, F3, c2, `m', `q');\\
$>$ A4:= sprem(F3, A3, c2, `m', `q');\\
$>$ A5:= sprem(A3, A4, c2, `m', `q');\\
$>$ sprem(A4, A5, c2, `m', `q');\\
$>$ factor(\%);\\}
\indent By factoring $R_1$:\\
{\tt $>$ factor(R1);\\
$>$ sprem(F3, c$1^\wedge$2 * c2 - 1, c2, `m', `q');\\
$>$ factor(\%);\\
$>$ sprem(F3, c2 - c1, c2, `m', `q');\\
$>$ factor(\%);\\
$>$ A6:= sprem(F3, c1 * c$2^\wedge$2 - 1, c2, `m', `q');\\
$>$ sprem(c1 * c$2^\wedge$2 - 1, A6, c2, `m', `q');\\
$>$ factor(\%);\\
$>$ A7:= sprem(R3, F3, c2, `m', `q');\\
$>$ A8:= sprem(F3, A7, c2, `m', `q');\\
$>$ A9:= sprem(A7, A8, c2, `m', `q');\\
$>$ sprem(A8, A9, c2, `m', `q');\\
$>$ factor(\%);\\}
\\
To solve $R_1=R_2=c^2_1 +c_1 +1=0$: \\
\indent We first solve $R_1=c^2_1 + c_1 + 1=0$.\\
{\tt $>$ B1:= sprem(R1, c$1^\wedge$2 + c1 + 1, c1, `m', `q');\\
$>$ sprem(c$1^\wedge$2 + c1 + 1, B1, c1, `m', `q');\\
$>$ factor(\%);\\}
\indent Then we solve $R_2=c^2_1 +c_1 +1=0$.\\
{\tt $>$ B2:= sprem(R2, c$1^\wedge$2 + c1 + 1, c1, `m', `q');\\
$>$ sprem(c$1^\wedge$2 + c1 + 1, B2, c1, `m', `q');\\
$>$ factor(\%);\\}
\\
To solve $R_1=R_2=c^2_1 - c_1 + 1 = 0$: \\
\indent We first solve $R_1=c^2_1 - c_1 +1=0$.\\
{\tt $>$ B3:= sprem(R1, c$1^\wedge$2 - c1 + 1, c1, `m', `q');\\
$>$ sprem(c$1^\wedge$2 - c1 + 1, B3, c1, `m', `q');\\
$>$ factor(\%); \\}
\indent Then we solve $R_2=c^2_1 - c_1 +1=0$.\\
{\tt $>$ B4:= sprem(R2, c$1^\wedge$2 - c1 + 1, c1, `m', `q');\\
$>$ sprem(c$1^\wedge$2 - c1 + 1, B4, c1, `m', `q');\\
$>$ factor(\%);\\}

\section{References}
\noindent $[1]$ J. Adriaenssens and L. Le Bruyn,
\textit{Non-commutative Covers and the Modular Group},
\textsf{math.RA/0307139} (2003). \\ \\
$[2]$ E. Letzter, \textit{Commutator Hopf Subalgebras and
Irreducible Representations}, ArXiv preprint math.RA/0507253.
\\ \\
$[3]$ M. Newman, \textit{The Structure of some Subgroups of the
Modular Group},  Illinois J. Math, \textbf{6} (1962) 480-487.
\\ \\
$[4]$ I. Tuba and H. Wenzl, \textit{Representations of the Braid
Group $B_3$ and of $SL(2,Z)$}, Pacific Journal of Mathematics,
\textbf{197} (2001) 491-510.
\\ \\
\textsc{Department of Mathematics, University of Wisconsin at Madison, Madison, Wisconsin, 53706} \\
\textit{E-mail address:} moran@math.wisc.edu \\ \\
\textsc{450 Memorial Drive, Cambridge, Massachusetts, 02139} \\
\textit{E-mail address:} matt\_tbo@mit.edu
\end{document}